\documentclass[10pt,a4paper]{article}
\usepackage[in]{newfullpage}
\usepackage{amsfonts}
\usepackage{amsthm}
\usepackage{amsmath}
\usepackage{amssymb}
\usepackage[all]{xy}
\usepackage{epsf}
\usepackage{epsfig}
\allowdisplaybreaks
\numberwithin{figure}{section}
\numberwithin{equation}{section}
\newtheorem{theor}{Theorem}[section]
\newtheorem{thm}[theor]{Theorem}
\newtheorem{lemma}[theor]{Lemma}

\newtheorem{dfn}[theor]{Definition}%[section]
\newtheorem{eg}[theor]{Example}%[section]
\newenvironment{example}{\begin{eg}\upshape}{\end{eg}}
\newtheorem*{note}{Note}
\newtheorem*{remark}{Remark}

\newcommand{\C}{{\mathbb C}}

\newcommand{\Z}{{\mathbb Z}}
\newcommand{\N}{{\mathbb N}}

\newcommand{\cB}{{\cal B}}
\newcommand{\cC}{{\cal C}}

\newcommand{\isom}{\cong}
\newcommand{\D}[1]{\ensuremath{D^{(#1)}}}
\newcommand{\T}[1]{\ensuremath{\mathrm{T}(#1)}}
\newcommand{\TT}[1]{\ensuremath{\widetilde{\mathrm{T}}(#1)}}
\newcommand{\set}[2]{\ensuremath{\{#1|#2\}}}
\newcommand{\Hom}{\ensuremath{\mathrm{Hom}}}
\newcommand{\Fun}{\ensuremath{\mathit{Fun}}}

\newcommand{\Root}{\ensuremath{\mathit{root}}}

\newcommand{\Multi}{\ensuremath{\mathit{Multi}}}
\newcommand{\Sing}{\ensuremath{\mathrm{Sing}}}
\newcommand{\Multibin}{\ensuremath{\mathit{Multi}}^K_{\epsfig{file=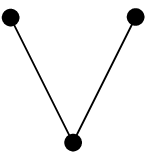,height=2mm}}}

\newcommand{\Partial}[1]{\ensuremath{\frac{\partial}{\partial #1}}}

\newcommand{\fhat}{\widehat{f}}
\newcommand{\ghat}{\widehat{g}}
\newcommand{\h}{\mbox{-}}

\newcommand{\tens}{\mbox{$\otimes$}}
\newcommand{\vac}{\mbox{$|0\rangle$}}

\newcommand{\Hn}{\ensuremath{H^{\otimes n}}}
\renewcommand{\H}[1]{\ensuremath{H^{\otimes #1}}}
\newcommand{\G}[1]{\ensuremath{G^{\otimes #1}}}

\newcommand{\preup}[2]{\ensuremath{\mbox{}^{#1}#2}}

\newcommand{\bottom}{\bot}

\newcommand{\acts}{\mbox{$\cdot$}}
\newtheorem{claim}[theor]{Claim}

\newcommand{\ps}[1]{[\![#1]\!]}

\newcommand{\bemph}[1]{\emph{\textbf{#1}}}
\newcommand{\defn}[1]{\textbf{#1}}

\CompileMatrices

\newcommand{\revar}[2]{\xymatrix{*=0{} \ar@<2pt>[r]^{#1} & *=0{}
\ar@<2pt>[l]^{#2}}} 
\newcommand{\revmapsto}[2]{\xymatrix{*=0{} \ar@{|->}@<3pt>[r]^{#1} & *=0{}
\ar@{|->}@<3pt>[l]^{#2}}} 

\newcommand{\genericblank}{\vcenter{
\xymatrix{ & & \ldots &  \\
	   & & *=0{} \ar@{-}[llu] \ar@{-}[lu] \ar@{-}[ru]}}}

\newcommand{\generic}[3]{\vcenter{
\xymatrix{{#1}_1&{#1}_2 &  \ldots	& {#1}_{#3}  \\
		& 	& *=0{} \ar@{-}[llu]^{{#2}_1}
			  \ar@{-}[lu]_{{#2}_2}\ar@{-}[ru]_{{#2}_{#3}}}}}

\newcommand{\genericnolabel}[3]{\vcenter{
\xymatrix{{#1}_1&{#1}_2 &  \ldots	& {#1}_{#2}  \\
		& 	& *=0{\txt{\\ \\{$#3$}}} \ar@{-}[llu] 
			  \ar@{-}[lu]\ar@{-}[ru]}}}

\newcommand{\flatmodule}[4]{\vcenter{
\xymatrix{{#1}_1&{#1}_2 &  \ldots	& {#1}_{#3}  \\
		& 	& *=0{\txt{\\ \\{$#4$}}} \ar@{-}[llu]^{{#2}_1}
			  \ar@{-}[lu]_{{#2}_2}\ar@{-}[ru]_{{#2}_{#3}}}}}

\newcommand{\flattwoleafed}[5]{\vcenter{
\xymatrix@C=10pt@R=15pt{{#1} &  & {#2}  \\
	   & *=0{\txt{\\ \\{$#3$}}} \ar@{-}[lu]^{#4} \ar@{-}[ru]_{#5}}}}

\newcommand{\flatthreeleafed}[4]{\vcenter{
\xymatrix@C=10pt@R=15pt{{#1} &{#2} & {#3} \\
	   & *=0{\txt{\\ \\{$#4$}}} \ar@{-}[lu] \ar@{-}[u] \ar@{-}[ru]}}}

\newcommand{\flatthreeleafedlabelled}[7]{\vcenter{
\xymatrix@C=10pt@R=15pt{{#1} &{#2} & {#3} \\
	   & *=0{\txt{\\ \\{$#4$}}} \ar@{-}[lu]^{#5} \ar@{-}[u]_{#6} \ar@{-}[ru]_{#7}}}}

\newcommand{\threeleafedleft}[8]{\vcenter{
\xymatrix@C=10pt@R=15pt{ {#1}&                               &{#2} & &  \\
	       & *=0{} \ar@{-}[lu]^{#5} \ar@{-}[ru]_{#6} &     & {#3}&  \\
               & & *=0{\txt{\\ \\{$#4$}}} \ar@{-}[lu]^{#7} \ar@{-}[ru]_{#8}}}}

\newcommand{\threeleafedright}[8]{\vcenter{
\xymatrix@C=10pt@R=15pt{     &                               &{#2} & &{#3}&  \\
	   &{#1}& & *=0{} \ar@{-}[lu]^{#5} \ar@{-}[ru]_{#6} &  \\
               & & *=0{\txt{\\ \\{$#4$}}} \ar@{-}[lu]^{#7} \ar@{-}[ru]_{#8}}}}

\newcommand{\fourleafdoubleA}{\vcenter{
\xymatrix@C=10pt@R=20pt{ V &  & V & & V & & V  \\
	     & *=0{} \ar@{-}[lu]^{x} \ar@{-}[ru] 
                &   & &   & *=0{} \ar@{-}[lu]^{y} \ar@{-}[ru] &  \\
             &  &   & *=0{\txt{\\ \\$V$}} \ar@{-}[llu]^{z} \ar@{-}[rru]}}}

\newcommand{\fourleafdoubleB}[5]{\vcenter{
\xymatrix@C=10pt@R=20pt{ {#1} &  & {#2} & & {#3} & & {#4}  \\
	     & *=0{} \ar@{-}[lu] \ar@{-}[ru] 
                &   & &   & *=0{} \ar@{-}[lu] \ar@{-}[ru] &  \\
             &  &   & *=0{\txt{\\ \\{$#5$}}} \ar@{-}[llu] \ar@{-}[rru]}}}

\newcommand{\fourleafdoubleC}[5]{\vcenter{
\xymatrix@C=10pt@R=20pt{ {#1}_1 &  & {#1}_2 & & {#1}_3 & & {#1}_4  \\
	     & *=0{} \ar@{-}[lu]^{{#3}_1} \ar@{-}[ru]_{{#3}_2}
                &   & &   & *=0{} \ar@{-}[lu]^{{#4}_1} \ar@{-}[ru]_{{#4}_2} &  \\
             &  &   & *=0{\txt{\\ \\{$#2$}}} \ar@{-}[llu]^{{#5}_1} \ar@{-}[rru]_{{#5}_2}}}}

\newcommand{\twoleavedrighttwo}[6]{\vcenter{
\xymatrix@C=10pt@R=15pt{    & & {#2}& &  \\
	                {#1}& & *=0{} \ar@{-}[u]^{#6} &  \\
                            & *=0{\txt{\\ \\{$#3$}}} \ar@{-}[lu]^{#4} \ar@{-}[ru]_{#5}}}}

\newcommand{\twoleafedtail}[6]{\vcenter{
\xymatrix@C=10pt@R=15pt{ {#1}&                                         &{#2}  \\
	                     & *=0{} \ar@{-}[lu]^{#4} \ar@{-}[ru]_{#5} &       \\
                             & *=0{\txt{\\ \\{$#3$}}} \ar@{-}[u]^{#6} }}}

\newcommand{\flatthreeleafedleft}[4]{\vcenter{
\xymatrix@C=10pt@R=15pt{ {#1}&                               &{#2} &{#3}&  \\
	       & *=0{} \ar@{-}[lu] \ar@{-}[ru] &     & *=0{} \ar@{-}[u]&  \\
               & & *=0{\txt{\\ \\$#4$}} \ar@{-}[lu] \ar@{-}[ru]}}}

\newcommand{\flatthreeleafeddown}{\vcenter{
\xymatrix@C=10pt@R=15pt{ & &  \\
	   & *=0{} \ar@{-}[lu] \ar@{-}[u] \ar@{-}[ru] \\
           & *=0{} \ar@{-}[u]  }}}

\newcommand{\genericA}[3]{\vcenter{
\xymatrix@C=10pt@R=15pt{{#1}_1&{#1}_2 &  \ldots	& {#1}_{#3}  \\
		& 	& *=0{\txt{\\ \\$#2$}} \ar@{-}[llu]
			  \ar@{-}[lu]\ar@{-}[ru]}}}

\newcommand{\generictailA}[3]{\vcenter{
\xymatrix@C=10pt@R=15pt{
{#1}_1	&{#1}_2	&\ldots	& {#1}_{#3}  \\
	&	&*=0{} \ar@{-}[llu]
			  \ar@{-}[lu] \ar@{-}[ru] \\
	&	&*=0{\txt{\\ \\$#2$}} \ar@{-}[u]}}}

\newcommand{\generictopA}[3]{\vcenter{
\xymatrix@C=10pt@R=15pt{
	&	&	&{#1}_k			&	&\\
{#1}_1	&{#1}_2	&\ldots	&*=0{} \ar@{-}[u]	&\ldots	&{#1}_{#3}  \\
	&	& *=0{\txt{\\ \\$#2$}} \ar@{-}[llu]
			  \ar@{-}[lu] \ar@{-}[ru] \ar@{-}[rrru]}}}

\newcommand{\genericcompositionB}{\vcenter{
\xymatrix@C=10pt@R=15pt{
A_{1}	& A_{2}	&\ldots	& A_{m}&\\
	&	& *=0{} \ar@{-}[llu] \ar@{-}[lu]\ar@{-}[ru]
			&A_{m+1}	&\ldots & A_n\\ 
	&		&&&*=0{\txt{\\ \\$B$}} \ar@{-}[llu]\ar@{-}[lu]
			 \ar@{-}[ru]  }}}

\newcommand{\genericcompositionA}{\vcenter{
\xymatrix@C=10pt@R=15pt{&& A_{i_1}& A_{i_2}	&\ldots	& A_{i_m}&\\
	  A_1&A_2&	&\ldots	& *=0{} \ar@{-}[llu]
			  \ar@{-}[lu]\ar@{-}[ru]&\ldots &A_n\\
	  &&	&	&*=0{\txt{\\ \\$B$}} \ar@{-}[llllu]
			  \ar@{-}[lllu] \ar@{-}[u]\ar@{-}[rru]  }}}

\newcommand{\flatthreeleafedleftsides}[6]{\vcenter{
\xymatrix{&                    & & &  \\
	  & *=0{} \ar@{-}[lu]^{#1} \ar@{-}[ru]_{#2} & & *=0{} \ar@{-}[u]_{#3}&  \\
          & & *=0{\txt{\\ \\{$#6$}\\}} \ar@{-}[lu]^{#4} \ar@{-}[ru]_{#5}}}}

\newcommand{\flatthreeleafedleftsidestop}[6]{\vcenter{
\xymatrix{A_1&                    &A_2 &A_3 &  \\
	  & *=0{} \ar@{-}[lu]^{#1} \ar@{-}[ru]_{#2} & & *=0{} \ar@{-}[u]_{#3}&  \\
          & & *=0{\txt{\\ \\{$#6$}\\}} \ar@{-}[lu]^{#4} \ar@{-}[ru]_{#5}}}}

\newcommand{\genericcomposition}{\vcenter{
\xymatrix{&&	&	&\ldots	& &\\
	  &&	&\ldots	& *=0{} \ar@{-}[llu]^{{x}_1}
			  \ar@{-}[lu]_{{x}_2}\ar@{-}[ru]_{{x}_{n}}&\ldots &\\
	  &&	&	&*=0{} \ar@{-}[llllu]^{{y}_1}
			  \ar@{-}[lllu]_{{y}_2} \ar@{-}[u]_{y_k}
\ar@{-}[rru]_{{y}_{m}}  }}}

\newcommand{\generictail}[3]{\vcenter{
\xymatrix{
{#1}_1	&{#1}_2	&\ldots	& {#1}_{#3}  \\
	&	&*=0{} \ar@{-}[llu]^{{#2}_1-y}
			  \ar@{-}[lu]_{{#2}_2-y} \ar@{-}[ru]_{{#2}_{#3}-y} \\
	&	&*=0{} \ar@{-}[u]^{y}}}}

\newcommand{\oldgenerictail}[3]{\vcenter{
\xymatrix{
{#1}_1	&{#1}_2	&\ldots	&{#1}_k	&\ldots	& {#1}_{#3}  \\
	&	& 	&*=0{} \ar@{-}[lllu]^{{#2}_1-{#2}_k}
			  \ar@{-}[llu]_{{#2}_2-{#2}_k} \ar@{-}[u]_{0}
\ar@{-}[rru]_{{#2}_{#3}-{#2}_k} \\
	&	&	&*=0{} \ar@{-}[u]^{{#2}_k}}}}

\newcommand{\generictop}[3]{\vcenter{
\xymatrix{
	&	&	&{#1}_k			&	&\\
{#1}_1	&{#1}_2	&\ldots	&*=0{} \ar@{-}[u]^{y}	&\ldots	&{#1}_{#3}  \\
	&	& *=0{} \ar@{-}[llu]^{{#2}_1}
			  \ar@{-}[lu]_{{#2}_2} \ar@{-}[ru]^{{#2}_k-y}
\ar@{-}[rrru]_{{#2}_{#3}}}}}

\newcommand{\genericflattop}[3]{\vcenter{
\xymatrix{
{#1}_1	&{#1}_2	&\ldots	&{#1}_k			&\ldots	&{#1}_{#3}  \\
	&	& *=0{} \ar@{-}[llu]^{{#2}_1}
			  \ar@{-}[lu]_{{#2}_2} \ar@{-}[ru]^{{#2}_k}
\ar@{-}[rrru]_{{#2}_{#3}}}}}

\newcommand{\threenleaf}{
\xymatrix{ 
a_1	& a_2 	& a_3 	& a_4 	& a_5	&  \ldots	& a_n  \\
	&*=0{}\ar@{-}[lu]_{y_1} \ar@{-}[u]_{y_2}
		&	& *=0{}\ar@{-}[lu]_{y_3}\ar@{-}[u]_{y_4}\ar@{-}[ru]_{y_5}
				&\ldots	& *=0{}\ar@{-}[ru]_{y_n} \\
	&	&	&	&*=0{}\ar@{-}[lllu]^{x_1} \ar@{-}[lu]_{x_2}
				\ar@{-}[ru]_{x_n}  
} }

\title{Equivalence of Borcherds $G$\h Vertex Algebras and Axiomatic Vertex
Algebras} 
\author{Craig T. Snydal\thanks{ctsnydal@dpmms.cam.ac.uk} }
\date{20 April 1999}

\begin{document}
\maketitle
\begin{abstract}

In this paper we build an abstract description of vertex algebras from their
basic axioms.  Starting with Borcherds' notion of a vertex group, we
naturally construct a family of multilinear singular maps parameterised by
trees.  These singular maps are defined in a way which focusses on the
relations of singularities to their inputs.  In particular we show that this
description of a vertex algebra allows us to present generalised notions of
rationality, commutativity and associativity as natural consequences of the
definition.  Finally, we show that for a certain choice of vertex group,
axiomatic vertex algebras correspond bijectively to algebras in the relaxed
multilinear category of representations of a vertex group.
\end{abstract}
\tableofcontents
\bigskip

\section{Introduction}

\subsection{Motivation} \label{ss:motivation}
The theory of (non-supersymmetric) axiomatic vertex algebras has as its data,
a complex vector space, $V$, together with a \defn{vertex operator}:
\[Y(\cdot,x)\cdot:V\tens V \longrightarrow V\ps{x}[x^{-1}],\]
as well as a distinguished vector $\vac$ and an automorphism $T:V\rightarrow
V$.  So for vectors $a,b \in V$, the vertex operator can be written as
\[Y(a,x)b = x^{-k} \sum_{i \geq 0} c_i x^i, \]
where $c_i \in V$ and $k \geq 0$.  We would like to describe carefully the
products of vertex operators, and to give a precise description of what types
of maps arise from such products.

We represent the collection of all such vertex operators by the labelled
tree
\[\flattwoleafed{V}{V}{V}{x}{}\]
By extending the domain of the vertex operator pointwise to
$V\ps{x}[x^{-1}]$, we have two ways to take the product of a pair of vertex
operators, giving maps:
\begin{eqnarray}
Y(\cdot,x)Y(\cdot,y)\cdot:&V\tens V\tens V \longrightarrow&
V\ps{x}[x^{-1}]\ps{y}[y^{-1}] \label{e:wrong_comp1}\\
Y(Y(\cdot,y)\cdot,x)\cdot:&V\tens V\tens V \longrightarrow&
V\ps{x}[x^{-1}]\ps{y}[y^{-1}].\label{e:wrong_comp2}
\end{eqnarray}
Immediately we notice that the collection of all maps from $V^{\otimes 3}$ to
$V\ps{x}[x^{-1}]\ps{y}[y^{-1}]$ contains many maps which can not be realised
as the composite of two vertex operators.  For example, the singularity in
the variable $y$ only depends on two of the copies of $V$ in $V^{\otimes 3}$.
Composing two copies of the labelled tree, we can represent the collection of
all maps arising as the composite of two vertex operators by
\[\threeleafedleft{V}{V}{V}{V}{y}{}{x}{}
\threeleafedright{V}{V}{V}{V}{y}{}{x}{} \] 
This labeled tree notation makes clear this dependency of singularities on
inputs.  Looking more closely at the product, $Y(\cdot,x)Y(\cdot,y)\cdot$, we
see that a more accurate description of the space of maps defined by this
product is given by the collection,
\[\Hom\Bigl(V\tens V,\Hom(V, V\ps{x}[x^{-1}])\ps{y}[y^{-1}]\Bigr).\]
We can repeat this process, so that the collection of $n$\h fold products of
vertex operators, $Y(\cdot,x_n)\cdots Y(\cdot,x_1)\cdot$ can be described
exactly as the product,
\[\Hom\biggl(V\tens V,\Hom\Bigl(V, \cdots \Hom(V, V\ps{x_n}[x_n^{-1}])\cdots
\Bigr)\ps{x_1}[x_1^{-1}]\biggr).\] 
Similarly, this space describes the other type of $n$\h fold product of
vertex operators, 
\[Y\Bigl(Y\bigl(\cdots Y(Y(\cdot,x_1)\cdot,x_{2})\cdots\bigr)\cdot,x_n\Bigr)\cdot.\]
The difficulty arises when we begin to consider the space of products of
vertex operators which contain both types of composition.  Take for example
the composite,
\[Y(Y(\cdot,x)\cdot,z)Y(\cdot,y)\cdot.\]
This product can be achieved either by composing three vertex operators,
$Y(\cdot, x)\cdot$, $Y(\cdot, y)\cdot$, and $Y(\cdot, z)\cdot$, or
by appropriately composing a vertex operator $Y(\cdot, y)\cdot$ with an
element of
\[\Hom(V\tens V,\Hom(V, V\ps{z}[z^{-1}])\ps{x}[x^{-1}])=
\threeleafedleft{V}{V}{V}{V}{x}{}{z}{}\] 
or even by composing $Y(\cdot, x)\cdot$ in the other way with an element of 
\[\Hom(V\tens V,\Hom(V, V\ps{z}[z^{-1}])\ps{y}[y^{-1}])=
\threeleafedright{V}{V}{V}{V}{y}{}{z}{}\] 
Either way, we would like to understand the space of all maps which can be
realized as any of the composites given above.  Using the tree
notation, we shall denote the space of all such maps by
\[\fourleafdoubleA \]
but this space can not be described as explicitly as the previous spaces of
products of vertex operators. 

In addition to the desire to simply describe these spaces of products of
vertex operators, we would like to be able to relate them in such a way as to
take account of the axioms of the vertex algebra.  The usual way of working
with products of vertex operators is to consider them inside very large
spaces of formal distributions.  While this has the benefit of leaving plenty
of room to manipulate power series, it tends to obscure the important
features of the vertex operators.  In this paper we shall take the opposite
tack, opting instead to use this minimal description of the spaces in which
these products of vertex operators live.  In this way, we will be able to
make clear the essential features of composition of vertex operators.

\subsection{Outline}
The idea behind this paper is the extraction the essential features of
axiomatic vertex algebras, with the ultimate goal of showing how they can be
seen to arise naturally in a certain context.  By demonstrating this
naturality, the generalisation of vertex algebras to more general settings
becomes simply a matter of choosing the correct categories in which to work,
with the long term goal of applying them to higher dimensional field
theories.

We shall approach the problem of giving an abstract account of vertex
algebras by first considering carefully the vertex operator.  After reviewing
the definition of a vertex algebra, we define a vertex group which ties
together the power series and the infinitesimal translation operator.
Temporarily setting aside the singularities of a vertex operator, we show how
holomorphic vertex algebras arise naturally when considering representations
of the vertex group.

Next we look at the role of singularities in the definition of a vertex
operator.  We show how to use the elementary vertex structure of our vertex
group to define a space of singular functions, and after considering
composition of these singular functions we define spaces of singular
functions parameterised by binary trees.  Our definitions, while motivated by
the Borcherds definition in \cite{bor}, are different because they emphasise
the relation of singularities on inputs.  For the classical vertex group, we
demonstrate the correspondence between a binary singular functions and vertex
operators.

In the next section we consider the axioms for our vertex algebra.  After
discussing the types of multi maps which arise for one leafed trees, we see
that the vacuum axioms lead us to consider maps associated to trees with zero
leaves.  We show how the locality condition for vertex operators suggests the
notion of singular maps associated to a three leafed tree with no internal
vertices.  We finish this section by defining singular maps associated to an
arbritrary tree, and that natural maps between $n$ leafed trees give rise to
natural maps between corresponding singular functions.

We finish this paper by describing the general categorical structure
satisfied by these singular maps.  They form a relaxed multilinear category,
and we show that a vertex algebra is just an algebra in this category.  

For clarity, our presentation does not give the supersymmetric details, but
the generalisation can be made easily.

\subsection{Notation}
Throughout this paper we will often refer to modules of polynomials, power
series, and Laurent series.  Given any $R$\h module, $B$ (where $R$ is a
commutative ring with 1), our convention will be to denote the collection of
power series with coefficients in $B$ by $B\ps{x}$, and the collection of
Laurent series by $B\ps{x}[x^{-1}]$.  The notation can be combined to form
larger modules such as $B\ps{x}[x^{-1}]\ps{y}$, the collection of power
series in the variable $y$, whose coefficients are Laurent series in the
variable $x$.  Note that $B\ps{x}[x^{-1}]\ps{y} \supsetneq B\ps{x,y}[x^{-1}]$
because, for example, $\sum_{i \geq 0} x^{-i} y^i$ is contained in the first
module but not the second.

We also adopt Sweedler's notation \cite{sweedler} to denote comultiplication
by $\Delta(h) = \sum_{(h)} h_{(1)} \tens h_{(2)}$ for any coalgebra element
$h$.  

\subsection*{Acknowledgements}
I am very grateful to Martin Hyland for his time and patience with this
project.  I would also like to thank Richard Borcherds for frequently
fielding my detailed questions about his idea of a vertex algebra.  And, a
special thanks goes to Tom Leinster for the many discussions related to the
completion of this paper, and especially for taking the time to show me his
definition of the category of trees.  Also a great thanks to the creators of
\Xy-pic who made the tree diagrams possible.

\section{Classical Vertex Group and its Representations} \label{s:classical_vg}

In the introduction we concentrated on examining the products of vertex
operators.  We purposely avoided the question of what types of homomorphisms
we were considering, because we were mainly concerned with demonstrating the
dependence of singularities on products of vertex operators.  In this section
we return to the vertex operator itself, examining it closely and paying
special attention to its interaction with the infinitesimal translation
operator.

Before beginning we recall the definition of a vertex algebra (following the
definition in \cite{kac}).  

\begin{dfn}\label{d:vetex_algebra}
A \defn{vertex algebra} consists of a complex vector space $V$ (the
\defn{state space}) together with an endomorphism, $T:V \rightarrow V$ (the
\defn{infinitesimal translation operator}), a distinguished vector
denoted $\vac \in V$ (the \defn{vacuum vector}), and a linear map (the
\defn{vertex operator}): 
\[Y(\cdot, x)\cdot:V \tens V \rightarrow V\ps{x}[x^{-1}],\]
taking any $a\tens b \in V\tens V$ to $Y(a,x)b \in V\ps{x}$.  
The axioms for a vertex algebra say that 
\begin{description}
\item[Vacuum axioms:] For any $a \in V$, the vacuum satisfies:
\begin{eqnarray}
T\vac &=& 0 \label{e:vac_triv_axiom}\\
Y(\vac, x)a &=& a\label{e:vac_id_axiom} \\
Y(a,x)\vac|_{x=0} &=& a\label{e:vac_ps_axiom}.
\end{eqnarray}
\item[Translation covariance axiom:] The infinitesimal translation operator
interacts with a vertex operator according to:
\begin{equation}
[T,Y(a,x)] = \partial_x Y(a,x):V \rightarrow V\ps{x}[x^{-1}].\label{e:inv_axiom}
\end{equation}
\item[Locality axiom:] And, for any $a, b \in V$, there exists some $N \gg 0$
such that the following holds:
\begin{equation}\label{e:locality}
(x-y)^N[Y(a,x), Y(b, y)] = 0.
\end{equation}
\end{description}
\end{dfn}

If we consider the free group ring generated by the infinitesimal translation
operator, denoted $G$, then $V$ is a $G$\h module simply by virtue of the fact that
$T$ is an endomorphism of $V$.  From a simple vertex algebra calculation,
one can show that $\partial_x Y(a,x) = Y(Ta,x)$ (see \cite[Prop. 4.8]{kac}).
So, from the translation covariance axiom, we see that given any $b \in V$,
we have an equality of Laurent series, 
\[Y(a,x)Tb + Y(Ta,x)b = TY(a,x)b.\]
Now, considering the vertex operator as a map from $V\tens V$ to
$V\ps{x}[x^{-1}]$, this presentation of the translation covariance axiom just
says that a vertex operator is a $G$\h invariant map, where $G$ acts on
$V\tens V$ by
\begin{eqnarray*}
G \tens (V \tens V) &\longrightarrow &V \tens V\\
T\tens a \tens b &\longmapsto &Ta \tens b + a \tens Tb.
\end{eqnarray*}

From this account, it is clear that a general vertex algebra is going to be
some type of module.  In the next section we will define a vertex group,
which will be the algebra over which we will want to work.  We show that $G$
is a nontrivial example of a vertex group.  We then move on to consider
modules over a vertex group, and modules over $G$ in particular.  We finish
this section by considering ways of expressing maps between $G$\h modules and
we show how holomorphic vertex algebras arise naturally from this presentation.

\subsection{Definition of a Vertex Group}
This section reviews some important definitions from \cite{bor}.  In
particular, we define an elementary vertex structure on a cocommutative Hopf
algebra which we use to introduce the notion of a vertex group as a Hopf
algebra together with a specific choice of elementary vertex structure.

Recall that a Hopf algebra is a module $H$ over a commutative ring $R$ (with
unit) that has both the structure of an algebra and a coalgebra.  It also
possesses an antipode map $S:H \rightarrow H$ which is both a map of algebras
and a map of coalgebras, and serves to connect the algebra and coalgebra
structures.

\begin{dfn}
Let $H$ be a Hopf algebra over a commutative ring $R$, and let $H^* =
\Hom_R(H,R)$ be the collection of $R$\h linear maps from $H$ to $R$.  We say
that an $R$\h module, $K$, gives an \defn{elementary vertex structure} on $H$
if it is an associative algebra over the algebra $H^*$ satisfying the
following properties:
\begin{enumerate}
\item \defn{(Closure under left and right translation)} $K$ is a 2\h sided
$H$\h module such that the product on $K$ is invariant under the left and
right actions of $H$.  By this we mean that for $h \in H$ and $k,l \in K$, we
have:
\begin{align*}
h\acts (kl)=\sum_{(h)}(h_{(1)}\acts k)(h_{(2)}\acts l) &&
(kl)\acts h=\sum_{(h)}(k\acts h_{(1)})(l\acts h_{(2)}).
\end{align*}
We also require that the unit map for the algebra, $\eta:H^* \rightarrow K$,
be a homomorphism of 2\h sided $H$\h modules.
\item \defn{(Closure under Inversion)} The antipode on $H$ gives rise to a
map $S^*$ on $H^*$ that can be extended to an $R$\h linear map on $K$.  By
abuse of notation, we will refer to this map as $S:K \rightarrow K$, and
because we are extending the the dual map it is not quite an $H^*$\h algebra
map, but instead satisfies $S(kl)=S(l)S(k)$ and $S(1)=1$ for $k,l,1 \in
K$.
\end{enumerate}
\end{dfn}

\begin{dfn}\label{d:K}
If $H$ is a cocommutative Hopf algebra and $K$ is an elementary vertex
structure on $H$ which is commutative and satisfies $S^2=1_K$, then we say
that we have an \defn{elementary vertex group}, $G$.  $H$ will be referred to
as the \defn{group ring} of the vertex group $G$ or the \defn{underlying Hopf
algebra}, and $K$ will be called the \defn{ring of singular functions on
$G$}.
\end{dfn}

\begin{note}
Since there are examples of vertex groups that are more general than this
definition allows, we have chosen the name elementary vertex groups.
Throughout the rest of this paper we will only be working with elementary
vertex groups, so for simplicity, we will refer to them simply as vertex groups.
\end{note}

\begin{dfn}
For any $H$\h module, $B$, we abuse terminology by calling $\Hom_R(H^n,B)$
the \defn{collection of nonsingular functions from $\G{n}$ to $B$}.  This
will also be written $\Hom_R(\G{n},B)$.
\end{dfn}

Notice that for any cocommutative Hopf algebra, $H$, letting $K=H^*$ gives
$H$ the structure of a trivial vertex group.  This is a 2\h sided $H$\h
module with left and right actions given as above for $g,h \in H$ and $f \in
H^*$,
\begin{align}\label{2sideaction}
(h\acts f)(g)=\sum_{(h)}h_{(1)}f(S(h_{(2)})g) && (f\acts h)(g)=f(hg).
\end{align}
Since we are working with the dual of a Hopf algebra this is closed under
left and right translation and under inversion.

\subsection{Classical Vertex Group}\label{ss:classicalvg}

The Hopf algebra which will prove most important for our later work is the
complex polynomial algebra in one variable, $H = \C[T]$.  For simplicity, we
shall denote powers of $T$ by $\D{i} = T^i/i!$.  Thus we have an associative
algebra with multiplication $\mu (\D{i} \tens \D{j}) = \binom{i+j}{i}
\D{i+j}$ for any $i,j \geq 0$ and unit $\D{0}$.  We shall denote
multiplication as $\D{i} \acts \D{j}$.

This module has the additional structure of a Hopf algebra from the following
three maps:
\begin{align*}
\text{Coalgebra Structure:} 	&& 
	&\begin{aligned}[t]\Delta: H &\rightarrow H \tens H\\\D{i} &\mapsto
	\sum_{p+q = i} \D{p}\tens\D{q}
	\end{aligned}
				&				
	&\begin{aligned}[t]  \epsilon:H &\rightarrow \C\\ \D{i} &\mapsto
	\begin{cases}
	1 & i=0 \\
	0 & \text{otherwise}
	\end{cases}
	\end{aligned}\\[5mm]
\text{Antipode:}	&&
	&\begin{aligned}[t]S: H &\rightarrow H.\\
	\D{i} &\mapsto (-1)^{i}\D{i}.
	\end{aligned}
\end{align*}
Note that in particular this structure is cocommutative (i.e., $(1 \tens
\Delta)\Delta = (\Delta \tens 1)\Delta$).  

It is easy to see that the linear dual of $H$, $H^* = \Hom_\C(H,\C)$, is the
ring of power series in one variable, $H^* = \C\ps{x}$.  The dual pairing is
given by letting $\D{1}$ act as differentiation on the $x$ variable evaluated
at $x=0$.  We can extend this action of $H$ on its dual to give an $H$\h
module structure on $H^*$ by defining the obvious action of derivation $\D{j}
\acts x^i = \binom{i}{j}x^{i-j}$, so that $H$ is just a group ring of
derivations on $\C\ps{x}$.  It is easy to see that more generally,
$\Hom_\C(\Hn,\C)$ can be identified with $\C\ps{x_1,\ldots,x_n}$, where a
map $f$ is paired with $\sum r_{i_1, \ldots, i_n} x_1^{i_1}, \ldots,
x_n^{i_n}$ when $f(\D{i_1} \tens \ldots \tens \D{i_n}) = r_{i_1, \ldots,
i_n}$.

If we choose our elementary vertex structure to be $K=\C\ps{x}[x^{-1}]$, the
space of Laurent series in one variable, we see immediately that it is an
associative algebra over $H^*$ by construction.  In other words, it is an
infinite dimensional vector space over $H^*$ spanned by $\set{x^{-i}}{i \geq
0}$, with multiplication and unit maps giving it the structure of an algebra.
The action of $H$ is the obvious extension of the action of derivation on
$H^*$.  So for any $j \geq 0$ and $i \in \Z$ we have $\D{j} \acts x^i =
\binom{i}{j}x^{i-j}$.  Recall that for negative $i$, the binomial coefficient
is given by
\[\binom{i}{j} = \frac{i(i-1) \cdots (i-j+1)}{j(j-1) \cdots 1} = (-1)^j
\binom{j-i-1}{j} .\]  
If $i=j$ then this gives back our dual pairing between $x^i$ and $\D{i}$.  We
extend the antipode to $K$ by defining $S(x^i)=(-1)^ix^i$ for all $i\in \Z$.
By the product rule for derivation, we also see that the action of $H$ on the
product of two elements in $K$ commutes with the multiplication by using the
diagonal map.  In other words, for any $j \geq 0$, $i_1,i_2 \in \Z$;
\begin{eqnarray*}
\D{j} \acts (x^{i_1} x^{i_2}) &=& \mu\bigl(\Delta(\D{j}) \acts (x^{i_1} \tens x^{i_2})\bigr)\\
&=& \sum_{p+q=j} \mu\bigl((\D{p} \acts x^{i_1}) \tens (\D{q} \acts x^{i_2})\bigr)\\
&=& \sum_{p+q=j} \binom{i_1}{p}\binom{i_2}{q} x^{i_1 + i_2 - p - q}\\
&=& \binom{i_1 + i_2}{j} x^{i_1 + i_2 - j}.
\end{eqnarray*} 
and so this vertex algebra is closed under left and right translation.  We
shall call this the classical vertex group since it will give rise to the
vertex algebras of \cite{kac}.  It shall be denoted $G$ and will be the most
important example for the sections that follow.

\subsection{Representations of the Classical Vertex Group}
We now return our attention to the case of a general vertex group, $G$, with
underlying Hopf algebra, $H$.  An $R$\h module, $B$, is a
\defn{representation} of the vertex group if it is a representation of the
the underlying Hopf algebra.  Hence the category of modules for a vertex
group, $G$, is exactly the category of $H$\h modules.  When we refer to the
action of $G$ on a module, we mean that its group ring is acting on the
module.  The reason for introducing the category of $G$\h modules is because
we will use the elementary vertex structure of the vertex group to provide
this category with some additional singular structure.  This additional
structure will make it into a relaxed multilinear category (see definition
\ref{d:relmultcat} below).  Before we introduce this additional structure we
look more closely at the underlying category of modules for a vertex
group.

\begin{note}
We will want to consider the collection of all representations of a vertex
group.  We shall not restrict our attention to the full subcategory of finite
dimensional representations because we want to allow ourselves the freedom to
work with representations possessing a singular structure freely generated as
an algebra.
\end{note}

We begin by pointing out that since the underlying Hopf algebra for any
vertex group is cocommutative, the category of $G$\h modules possesses a
symmetric tensor product, which is just the tensor product of objects
considered as $R$\h modules.  Given $G$\h modules $A$ and $B$, the action of
$G$ on $A \tens B$ is given by the diagonal map.  So we can think of $A \tens
B$ as having both a $G$\h action and a $(G \tens G)$\h action.  More
generally, given the tensor product of a collection of $G$\h modules, $A_1
\tens \cdots \tens A_n$, for every $1 \leq i \leq n$ there exist actions
of $\G{i}$ corresponding to the different ways of comultiplying $\G{i}
\rightarrow \G{n}$.

Next we notice that the collection of $R$\h linear maps between any two $G$\h
modules $A$ and $B$, denoted $\Hom_R(A,B)$, can be given the structure of a
module in a number of ways.  Given any such map, $f:A \rightarrow B$, we
first define an action of $G$ on $f$ for any $a \in A$ and $g \in
G$ by 
\begin{equation} \label{e:hom_action}
(g \acts f)(a) = \sum_{(g)}g_{(1)}\acts f (S(g_{(2)})\acts a).
\end{equation}
We say that $f$ is \defn{invariant} under the action of $G$ if $g \acts
f = \epsilon(g) f$ for all $g \in G$.  This is equivalent to saying $f(g
\acts a) = g \acts f(a)$, or in other words the $G$\h invariant ring maps are
exactly the $G$\h module morphisms.

There also exist two additional actions of $G$ on $\Hom_R(A,B)$; one where
$G$ acts on the domain, the other where it acts on the codomain.  The $G$\h
action on the domain of a map $f:A \rightarrow B$ is a right action on $f$,
while the $G$\h action on the codomain of $f$ is a left action.  In fact,
considering $\Hom_R(A_1 \tens \ldots \tens A_n,B)$ we see that in addition to
being a map of $G$\h modules, the domain possesses the additional structure
of a module for $\G{n}$.  So $\Hom_R(A_1 \tens \ldots \tens A_n,B)$ possesses
the structure of a (right) $\G{n}$\h module by action on the domain of $f$
(and similarly, a (left) $G$\h module for action on the codomain, $B$).  The
fact that $\Hom_R(A_1 \tens \ldots \tens A_n,B)$ can be a module for both $G$
and $\G{n}$ will be important when we look at the composition of singular
maps below (see section \ref{ss:sing_mult_maps}).

Returning to the particular case of the classical vertex group, we have shown
that $\Hom_\C(G,\C) \isom \C\ps{x}$.  We see immediately that for any $G$\h
module, $B$, the linear dual, $\Hom_\C(G,B)$, is isomorphic as a $G$\h module
to $B\ps{x}$, where $B\ps{x}$ denotes the collection of power series in $x$
with coefficients in $B$.  This isomorphism between maps linear maps and
power series can be extended to maps from $n$\h fold products of $G$ for any
$n \geq 0$.  In such a case we have $\Hom_\C(\G{n},B) \isom B \ps{x_1,
\ldots, x_n}$ (where we have tacitly included the case where $n=0$ and
$\Hom_\C(\G{0},B) = \Hom_\C(\C,B) \isom B$).

From the discussion above we know that $\G{n}$ has the structure of a $G$\h
module using the $n$\h fold diagonal map, so the collection of linear
maps $\Hom_\C(\G{n},B)$ also has the structure of a $G$\h module as described
above.  Under the isomorphism $\Hom_\C(\G{n},B) \isom B \ps{x_1, \ldots, x_n}$,
the action of $\D{i} \in G$ on any $b(x_1, \ldots, x_n) \in B \ps{x_1,
\ldots, x_n}$ is given by:
\begin{equation}\label{e:action}
\D{i} \acts b(x_1, \ldots, x_n) = \sum_{j_0 + \cdots + j_n = i}
(-1)^{j_0}\bigl(\D{j_0}_B + \D{j_1}_{x_1} \cdots + \D{j_n}_{x_n}\bigr)
b(x_1,\ldots,x_n).
\end{equation}
In the previous expression, $D_B$ denotes the action of $G$ on the
coefficients of the power series $b(x_1, \ldots, x_n)$, and where
$\D{j_k}_{x_k}$ denotes the action of differentiation on the variable $x_k$.
As before, we say that $b(x_1, \ldots, x_n)$ is invariant under the action of
$H$ (or $G$) if $\D{i} \acts b(x_1, \ldots, x_n)$ is zero.  This is the same
as requiring that $b(x_1, \ldots, x_n)$ satisfy the following for all $i \geq
0$:
\begin{equation}\label{e:invariant} 
\D{i}_B b(x_1, \ldots, x_n) = \sum_{j_1 + \cdots + j_n = i} \bigl(\D{j_1}_{x_1} +
\cdots + \D{j_n}_{x_n}\bigr) b(x_1, \ldots, x_n).
\end{equation}
The collection of all such invariant elements is denoted $\Hom_G(\G{n},B)$.

Working out equation \eqref{e:action} explicitly for the case of
$\Hom_\C(\G{2},B) \isom B \ps{x,y}$, we have
\begin{eqnarray*}
\D{1}\acts(\sum_{i,j \geq 0} b_{i,j} x^{i} y^{j}) &=&
\sum_{i,j \geq 0} (\D{1}\acts b_{i,j}) x^{i} y^{j} - 
(\D{1}_{x} + \D{1}_{y})\sum_{i,j \geq 0} b_{i,j} x^{i} y^{j} \\
&=& \sum_{i,j \geq 0} (\D{1}b_{i,j}) x^{i} y^{j} - 
\sum_{i,j \geq 0} i b_{i,j} x^{i-1} y^{j} - \sum_{i,j \geq 0} j b_{i,j} x^{i} y^{j-1}. 
\end{eqnarray*}
And from equation \eqref{e:invariant}, we see that the $G$\h invariant
elements of $B\ps{x,y}$ are exactly those which satisfy $\D{1} \acts b_{i,j}
= (i+1)b_{i+1,j} + (j+1)b_{i,j+1}$, or more generally
\begin{equation}
\D{k}\acts b_{i,j} = \sum_{p+q=k} \binom{p+i}{i} \binom{q+j}{j}b_{p+i,q+j}.
\end{equation}
For any $n$, this equation also has an ``inverted form'' in which
we can write any $b_{i_1,\ldots,i_n}$ as a sum of $b_{j_1,\ldots,j_n}$ where
one of the indices $(j_1,\ldots,j_n)$ is set to zero.  For $n=2$ we can write
$b_{i,j}$ as:
\begin{eqnarray}
b_{i,j} &=& \sum_{p+q=j} (-1)^p \binom{i+p}{i}\preup{\D{q}}{b_{i+p,0}}\label{eq:d_i0}
\\	&=& \sum_{p+q=i} (-1)^q \binom{j+q}{j}\preup{\D{p}}{b_{0,j+q}}.\label{eq:d_0i}
\end{eqnarray}
This ability to represent $G$\h invariant power series in various way will be
important for encapsulating duality for vertex algebras.

\subsection{Extended Representation of Morphisms - Holomorphic Vertex
Algebras}\label{ss:extended} 
Now that we have a characterisation of modules for a vertex group, we shall
look more closely at maps between them.  We would like to describe module
maps in a way which will allow us to naturally add the elementary vertex
structure.  This treatment will avoid a very general abstract treatment,
instead aiming to provide the flavour of theory.  The details can be found in
\cite{cts_valgdetails}.

Given any linear map between $G$\h modules, $f:A \rightarrow B$, we
define a unique $G$\h linear map $\fhat:A \rightarrow \Hom_R(G,B)$ by:
\[ \fhat(a)(g) = f(g \acts a). \]
This pairing between linear maps and $G$\h linear maps is actually an
isomorphism in the category of $G$\h modules, so we may easily pass from one
description to the other.  We shall call this an \defn{extended
representation} of $f$.  If the map $f$ was also $G$\h invariant, then the
extended representation of $f$ would be a $G$\h linear map from $A$ to $G$\h
invariant elements of $\Hom_R(G,B)$.  These extended representations can also
be used to reexpress linear maps of the form $f:A_1 \tens \cdots \tens A_n
\rightarrow B$ as $\G{n}$\h module maps $\fhat: A_1 \tens \cdots \tens A_n
\rightarrow \Hom_R(\G{n},B)$, where the action of $\G{n}$ on the domain, $A_1
\tens \cdots \tens A_n$, passes to an action on the domain of the maps
$\Hom_R(\G{n},B)$. 

\begin{example}
In the case of the classical vertex group, this means that linear maps
between $G$\h modules $f:A_1 \tens \cdots \tens A_n \rightarrow B$, are the
same as $\G{n}$\h linear maps
\[\fhat: A_1 \tens \cdots \tens A_n \rightarrow B\ps{x_1, \ldots, x_n}.\]
In particular, we have that the action of $G$ on any of the $A_i$ carries
over to differentiation of $x_i$ as
\[ \fhat(a_1\tens \cdots \tens \D{1}a_i \tens \cdots \tens a_n) =
\partial_{x_i} \fhat(a_1\tens \cdots \tens a_i \tens \cdots \tens a_n).\]
\end{example}

The composition of maps in this extended representation is slightly delicate.
Given ordinary linear maps of $G$\h modules, say $f:A \rightarrow B$ and $g:B
\rightarrow C$, then clearly these compose to give a map $g \circ f:A
\rightarrow C$.  Now $f$ and $g$ can be considered in their extended
representation, i.e., $G$\h linear maps $\fhat:A \rightarrow \Hom_R(G,B)$ and
$\ghat:B \rightarrow \Hom_R(G,C)$, which compose to give a $G$\h linear map
$\ghat \circ \fhat:A \rightarrow \Hom_R(\G{2}, C)$.  But the ordinary composite
$g \circ f$ has an extended representation $\widehat{g \circ f}:A \rightarrow
\Hom_R(G, C)$.  We would expect these to be related, but in fact, we need to
require that $f$ be $G$\h invariant in order to relate/reduce the extra
factor of $G$ in the codomain of $\ghat \circ \fhat$.  Intuitively, this
problem appears because the action of $G$ on $B$ appears in $\Hom_R(G,C)$,
but the composite $g \circ f$ ``hides'' that action.  So by taking $f$ to be
$G$ \h invariant we eliminate the action of $G$ on $f$ explicitly.

More generally, maps in the extended representation compose pointwise as
multilinear functions in the usual way, provided that all maps (except
possibly the bottom map) are $G$\h invariant.  

\begin{example}\label{ex:nonsing_comp}
Again we consider the case of the classical vertex group.  Let $A_i, B_j,
C_l$ and $D$ be $G$\h modules.  Given a $\G{n}$\h linear map $g:C_1 \tens
\cdots \tens C_n \rightarrow D\ps{z_1, \ldots, z_n}$, and $G$\h invariant
maps $f:A_1 \tens \cdots \tens A_p \rightarrow C_1\ps{x_1, \ldots, x_p}$,
$k:B_1 \tens \cdots \tens B_q \rightarrow C_2\ps{y_1, \ldots, y_q}$
($\G{p}$\h linear and $\G{q}$\h linear respectively), then the
composite is a $\G{(p+q+n-2)}$\h linear map
\[g \circ (f \tens k):A_1 \tens \cdots \tens A_p \tens B_1 \tens \cdots \tens
B_q \tens C_3 \tens \cdots \tens C_n \rightarrow D\ps{x_1, \ldots, x_p, y_1,
\ldots, y_q, z_1, \ldots, z_n}\]
which satisfies $\partial_{z_1} = \partial_{x_1}+ \cdots + \partial_{x_p}$
and $\partial_{z_2} = \partial_{y_1}+ \cdots + \partial_{y_q}$.  In other
words, the map $g \circ (f \tens k)$ factors through $D\ps{x_1+z_1, \ldots,
x_p+z_1, y_1+z_2, \ldots, y_q+z_2, z_3, \ldots, z_n}$.  Also, it makes sense
to consider compositions in a particular order, say $(g \circ k) \circ f$ or
$(g \circ f) \circ k$, and these are equal.
\end{example}

\begin{claim}\label{c:holo_equivalence}
Let {\boldmath $H \mbox{-} \mathit{Mod}$} be the category of representations
of the underlying Hopf algebra, $H$, of the classical vertex group. Then
category of holomorphic vertex algebra is isomorphic to the category of $H$\h
invariant commutative algebras in {\boldmath $H \mbox{-} \mathit{Mod}$}.
\end{claim}

\begin{proof}
In order to make sense of this statement, we will first review both the
definition of a holomorphic vertex algebra, and the definition of an algebra
in a category.  

\begin{dfn}
A \defn{holomorphic vertex algebra} is a vertex algebra without
singularities.  In other words, the vertex operator of a holomorphic vertex
algebra is a map:
\[Y(\cdot, x)\cdot:V \tens V \rightarrow V\ps{x}.\]
It follows that the locality axiom reduces to the statement that products of
vertex operators commute.

The collection of holomorphic vertex algebras is given the structure of a
category by defining a morphism in the category of holomorphic vertex
algebras to be a map of complex vector spaces taking vacuum to vacuum,
commuting with multiplication, and respecting the actions of the
infinitesimal translation operators.
\end{dfn}

\begin{dfn}\label{d:algebra_defn}
A \defn{commutative (associative) algebra} in a symmetric tensor category
consists of an object $A$, a multiplication map $\mu: A \tens A \rightarrow
A$ which is invariant under the symmetry action for the tensor product, and a
unit for the multiplication $\eta:I \rightarrow A$ satisfying the usual
axioms for associativity and unit (where $I$ is the unit for the tensor
product).
\end{dfn}

For the category of representations of the Hopf algebra $H$, the unit for
multiplication is $\C$. and the multiplication map $\mu$ is $H$\h invariant.
Given any such algebra in this category of representations, where the
multiplication map $\mu$ is $H$\h invariant, we form a holomorphic vertex
algebra as follows.  We begin by taking $T = \D{1}$ as the infinitesimal
translation operator, and take the vacuum vector to be $\eta(1)$.  Then by
considering $\mu$ in the extended representation, we have a map
\[ \widehat{\mu}:A\tens A \rightarrow A\ps{x,y},\]
which satisfies 
\begin{eqnarray}
\widehat{\mu}(Ta\tens b) &=& \partial_x \widehat{\mu}(a \tens b)\\
\widehat{\mu}(a\tens Tb) &=& \partial_y \widehat{\mu}(a \tens b)\\
T\widehat{\mu}(a \tens b)&=& \widehat{\mu}(Ta\tens b)+\widehat{\mu}(a\tens Tb).
\end{eqnarray}
From this we define a vertex operator on $A$ to be 
\begin{equation}\label{e:operator_from_alg}
Y(\cdot, x)\cdot = \widehat{\mu}(\cdot \tens \cdot)|_{y=0}:A \tens A
\rightarrow A\ps{x}.
\end{equation}

Checking that this satisfies the axioms for a holomorphic vertex algebra, we
see immediately that the vacuum axioms are satisfied because $H$ acts
trivially on $\C$, and because $\mu(\eta(1)\tens a) = a$ for all $a \in V$.
Writing out the translation covariance axiom, we use the $H$\h invariance of
$\widehat{\mu}$ to give:
\begin{eqnarray}
T\widehat{\mu}(a \tens \cdot)|_{y=0} - \widehat{\mu}(a \tens T\cdot)|_{y=0} &=&
\widehat{\mu}(Ta\tens \cdot)|_{y=0} +\widehat{\mu}(a\tens T\cdot)|_{y=0} -
\widehat{\mu}(a \tens T\cdot)|_{y=0} \\
&=& \widehat{\mu}(Ta\tens \cdot)|_{y=0} \\
&=& \partial_x \widehat{\mu}(a\tens \cdot)|_{y=0}.
\end{eqnarray}
And the locality axiom follows from the commutativity of $\widehat{\mu}$.
Notice that $\widehat{\mu}(a \tens b) = Y(a,x)Y(b,y)\vac$.

Similarly, given any vertex algebra, we can easily define an algebra in the
category of of representations of the Hopf algebra $H$ by taking $\eta(r) =
r\vac \in V$ for any $r \in R$, and $\mu(\cdot \tens \cdot) = Y(\cdot,
x)\cdot|_{x=0}$.  The unit axioms follow from the properties of the vacuum,
and associativity follows from locality.  The translation covariance axiom
says that the multiplication is $H$\h invariant.  And, locality axiom acting
on the vacuum says that this multiplication is commutative.

If $(V, Y(\cdot,x)\cdot, T, \vac)$ is a vertex algebra, and the corresponding
$H$\h invariant algebra is $(V, \mu, \eta)$, then from this algebra, we get
back the vertex algebra, $(V, Y^\prime(\cdot,z)\cdot, T, \vac)$.  If
$Y(a,x)b = \sum_{i \geq 0} c_i x^i$ for $a\tens b \in V\tens V$, then 
\[\widehat{\mu}(a\tens b) = \sum_{i, j \geq
0}\biggl[Y(\D{i}a,x)\D{j}b\biggr]_{x=0} z^i y^j \] 
so setting $y=0$, we have
\begin{eqnarray*}
\widehat{\mu}(a\tens b)|_{y=0} &=& \sum_{i \geq
0}\biggl[Y(\D{i}a,x)b\biggr]_{x=0} z^i \\ 
&=& \sum_{i \geq 0}\biggl[\sum_{j
\geq 0} \binom{j+i}{i}c_{j+i} x^j\biggr]_{x=0} z^i \\ 
&=& \sum_{i \geq 0}c_i z^i,
\end{eqnarray*}
and so $Y = Y^\prime$.  Similarly, starting with an $H$\h invariant algebra,
$(A,\mu, \eta)$, the process of creating a holomorphic vertex algebra and
then mapping back to an $H$\h invariant algebra clearly maps $(A,\mu, \eta)$
to itself.

Finally, a morphism in the category of holomorphic can be seen to correspond
exactly to a morphism of algebras in the category of $H$\h modules.  Because
this pairing of objects in each category is completely natural we have an
isomorphism of categories and so our claim is proved.
\end{proof}

\section{Singularities}

Up to this point, we have concentrated on the underlying Hopf algebra
structure of our vertex group.  We now use the elementary vertex structure of
the vertex group to add singularities to the extended representations of
morphisms of $G$\h modules, paying special attention to the case for the
classical vertex group.

It is in this chapter that this paper departs from Richard Borcherds' paper
\cite{bor}.  The ideas in this section were inspired by that paper, but the
collections of singular maps we shall define are much smaller than those
defined in his paper.  In fact, they will be defined in such a way as to make
their relationship to his singular maps clear.  

\subsection{Localisation}\label{localisation}
The following definitions illustrate how we will put together the space of
nonsingular functions from $\G{n}$ to $B$ with the ring of singular functions
on $G$ in order to arrive at a notion of singular functions from $\G{n}$ to
$B$.  These will also be referred to as \emph{singular functions of type
$K$}.

\begin{dfn}
For $1 \leq i < j \leq n$ and any (right) $\Hom_R(\G{n},R)=(\G{n})^*$\h
module, $M$, the \defn{localisation of $M$ at $(i,j)$} is defined to be $M
\tens_{G^*} K$, where $M$ is given the structure of an $G^*$\h module by the
$R$\h module dual of the map
\begin{equation}\label{e:f_ij}
\begin{aligned} 
\widehat{f}_{ij}&: && \ \G{n} \ \ \longrightarrow \ \ &&  \ \ G. \\
&g_1 \tens &&\cdots \tens g_n \mapsto &&g_iS(g_j) 
\end{aligned}
\end{equation}
\end{dfn}

We have already seen that the space of nonsingular functions from $\G{n}$ to
$B$ has the structure of both an $G$\h module and an $\G{n}$\h module.
Because of the coassociativity of $G$, it can also be considered as an
$(\G{n})^*$\h module where the action of any $\alpha \in (\G{n})^*$ on
$f:\G{n} \rightarrow B$ is
\begin{equation} \label{e:dual_act}
(\alpha\cdot f)(h_1 \tens \cdots \tens h_n) = \alpha(h_1 \tens \cdots \tens
h_n)\cdot f(h_1 \tens \cdots \tens h_n)
\end{equation}
for all $h_1 \tens \cdots \tens h_n \in \G{n}$.  With this action defined, it
makes sense to localise the space of nonsingular functions from $\G{n}$ to
$B$.

\begin{lemma}\label{l:sequentail_locs}
The localisation of $\Hom_R(\G{n},B)$ is also a $(\G{n})^*$\h module,
and sequential localisations commute.
\end{lemma}

\begin{proof}
This is clear from the action defined in equation \eqref{e:dual_act} and the
coassociativity of $G$.
\end{proof}

\begin{dfn}\label{d:fun}
The \defn{space of singular functions} from $\G{n}$ to $B$, denoted
$\Fun(\G{n},B)$, is defined to be the localisation of the space of nonsingular
functions from $\G{n}$ to $B$ at all $(i,j)$ for $1 \leq i < j \leq n$.  It is
the space 
\[\Hom_R(\G{n},B) {\bigotimes}_{\substack{f_{i,j}\\1 \leq i < j \leq n}} K \]
\end{dfn}

The space of singular functions from $\G{n}$ to $B$ is a module for a number
of different actions.  Firstly, from lemma \ref{l:sequentail_locs}, we have
that $\Fun(\G{n},B)$ is a $(\G{n})^*$\h module.  Next we know that there is a
left action of $G$ on $B$ which carries over to $\Fun(\G{n},B)$.  We also
know that there is a (left) diagonal action of $G$ on the domain of the
nonsingular maps $\Hom_R(\G{n},B)$.  An easy check shows that this diagonal
action carries over to a trivial action on $K$ via the maps $f_{i,j}$ defined
in equation \eqref{e:f_ij}.  Hence the standard action of $G$ on any singular
function is the same as the standard action of $G$ on the nonsingular
functions given by \eqref{e:hom_action}.  This gives us the following lemma:

\begin{lemma}
If we define the \defn{space of singular $G$\h invariant functions} from
$\G{n}$ to $B$, denoted $\Fun_G(\G{n},B)$, to be the localisation of
$\Hom_G(\G{n},B)$ at all $(i,j)$ for $1 \leq i < j \leq n$, then this
coincides with the $G$\h invariant elements of $\Fun(\G{n},B)$ under the
standard action of $G$.
\end{lemma}
Finally we know that the nonsingular maps, $\Hom(\G{n},B)$, possesses a (right)
action of $\G{n}$ on their domain.  This action can be extended to singular
maps by using comultiplication to act on each of the $\otimes_{f_{i,j}} K$
terms through that map $f_{i,j}$.

\begin{example}
If we are working with a trivial vertex group (i.e., $K = H^*$)
then the localisation of the nonsingular functions from $\G{n}$ to $B$ has no
effect, and we just have $\Fun(\G{n},B) = \Hom_R(\H{n},B)$.
\end{example}

\begin{example}\label{ex:G1}
If we let $G$ be the classical vertex group then we showed in section
\ref{ss:classicalvg} that the collection of nonsingular functions,
$\Hom_\C(G,B)$, corresponds to $B\ps{x}$.  Since this can not be localised, we
have $\Fun(G,B)= \Hom_\C(G,B) \isom B\ps{x}$.  The $G$\h invariant maps are
given by $\Fun_G(G,B)=\Hom_G(G,B) \isom B$.
\end{example}

\begin{example}\label{ex:G2}
Again letting $G$ be the classical vertex group, we have also shown that the
collection of nonsingular functions from $\G{2}$ to $B$ corresponds to
$B\ps{x_{1},x_{2}}$, which we consider as a module over $\C\ps{x_{1},x_{2}}$.
The ring of singular functions on $G$ is isomorphic to $\C\ps{z}[z^{-1}]$.  We
can only localise at $(1,2)$ so consider the dual of the map
$\widehat{f}_{12}:\H{2} \rightarrow H$.  For any $\D{a},\D{b} \in H$, this
gives $\widehat{f}_{12}(\D{a}\tens \D{b}) = \binom{a+b}{a}(-1)^b \D{a+b}$, so the
dual map $f_{12}:\C\ps{z} \rightarrow \C\ps{x_{1},x_{2}}$ takes $z^k \in H^*$
to
\begin{eqnarray*} 
f_{12}(z^k) & = & \sum_{p+q=k} \binom{p+q}{p} (-1)^q x_{1}^p x_{2}^q \\   
	    & = & (x_{1}-x_{2})^k.
\end{eqnarray*}
This map can be extended from $H^*$ to include negative values of $k$, in
which case we will have infinitely many nonzero summands for negative values
of $q$, and will still have $f_{12}(z^k)=(x_{1}-x_{2})^k$.  Therefore our space of
singular functions from $\G{2}$ to $B$ is given by
\begin{eqnarray*} 
\Fun(\G{2},B) & = & B\ps{x_{1},x_{2}} \tens_{H^*} \C\ps{z}[z^{-1}] \\
	&\isom& B\ps{x_{1},x_{2}}\ps{(x_{1}-x_{2})}[(x_{1}-x_{2})^{-1}] \\
	&=& B\ps{x_{1},x_{2}}[(x_{1}-x_{2})^{-1}].
\end{eqnarray*} 
The $G$\h invariant subcollection of maps in $\Fun(\G{2},B)$ is the quotient
of $B\ps{x_{1},x_{2}}[(x_{1}-x_{2})^{-1}]$ by the relation $\D{1}_B
-\Partial{x_{1}} - \Partial{x_{2}}$, where $\D{1}_B$ is the action of $\D{1}$
on the module $B$.  We thus have $\Fun_G(\G{2},B) \isom
B\ps{x_{1}-x_{2}}[(x_{1}-x_{2})^{-1}]$.
\end{example}

\begin{example}\label{ex:G3}
Continuing the previous example for the case of $\G{3}$, we know that the space
of nonsingular functions from $\G{3}$ to $B$ corresponds to
$B\ps{x_{1},x_{2},x_{3}}$ and we need to localise at the three ordered pairs
(i,j) for $1 \leq i < j \leq 3$.  The dual maps $f_{12}$, $f_{13}$, $f_{23}$,
behave exactly as in the previous example, and so we have for our space of
singular functions from $\G{3}$ to $B$ the collection
\begin{eqnarray*}
\Fun(\G{3},B) &=& \begin{aligned}[t]B\ps{x_{1},x_{2},x_{3}}\tens \C
  \ps{(x_{1}-x_{2})}&[(x_{1}-x_{2})^{-1}] \tens \C
  \ps{(x_{1}-x_{3})}[(x_{1}-x_{3})^{-1}] \\&\tens \C
  \ps{(x_{2}-x_{3})}[(x_{2}-x_{3})^{-1}]
  \end{aligned}\\
&\isom&B\ps{x_{1},x_{2},x_{3}}[(x_{1}-x_{2})^{-1},(x_{1}-x_{3})^{-1},
(x_{2}-x_{3})^{-1}]. 
\end{eqnarray*}
This can be extended in the obvious way to functions from $\G{n}$ to $B$
giving:   
\begin{equation}
\Fun(\G{n}, B) \isom B\ps{x_1, \ldots, x_n}[(x_i-x_j)^{-1}_{1 \leq i < j \leq
n}]. 
\end{equation}
\end{example}

\subsection{Singular Multilinear Maps}\label{ss:sing_mult_maps}
Now that we know how to use the elementary vertex structure of our vertex
groups to add singularities to the picture, we can easily see how to
generalise the extended representations of multilinear maps to define
singular multilinear maps.  But following the discussion given in section
\ref{ss:motivation}, we want to define these maps very delicately.  

Consider first the collection of linear maps between $G$\h modules $A$ and
$B$, $\Hom_R(A,B)$.  We saw in section \ref{ss:extended} that this collection
is isomorphic to the collection of $G$\h linear maps from $A$ to
$\Hom_R(G,B)$.  Replacing the collection of nonsingular maps $\Hom_R(G,B)$
with the collection of singular maps $\Fun(G,B)$ has no effect since they are
equal, so we have:

\begin{dfn}\label{d:sing_1and2}
The singular multilinear maps from $A$ to $B$ are just the
multilinear maps in the extended representation:
\[\Multi^K(A;B) = \Hom_{G}(A, \Hom(G,B)).\]

Repeating the same process for the collection of maps between $G$\h modules
$A_1\tens A_2$ and $B$, we have the following as our singular multi maps:
\[\Multi^K(A_1,A_2;B) = \Hom_{\G{2}}(A_1\tens A_2, \Fun(\G{2},B)).\]
\end{dfn}

\begin{lemma}\label{l:singular_inv}
The $G$\h invariant singular multilinear maps, $\Multi_G^K(A_1, A_2;B)$ are
the $\G{2}$\h invariant maps from $A_1\tens A_2$ to $\Fun_G(\G{2},B)$.
\end{lemma}

\begin{proof}
The $\G{2}$\h invariance of the singular multilinear maps allows the action
of $G$ on the domain, $A_1\tens A_2$, to carry over to an action on the
domain of $\Fun(\G{2},B)$ as described in section \ref{ss:extended}.  Hence
the action of $G$ on $A_1\tens A_2$ passes through to an action on $B$
exactly when $\Fun(\G{2},B)$ is $G$\h invariant.
\end{proof}

\begin{example}\label{ex:sing_1_hom}
Take $G$ to be the classical vertex group, and let $f$ be a singular
multilinear maps in $\Multi^K(A;B) = \Hom_{G}(A,B\ps{x}).$ Then for any $a
\in A$ we can write $f(a) = \sum_{i\geq 0} b_i x^i$ for $b_i \in B$.  If $f$
is $G$\h invariant, then $Tf(a) = \partial_xf(a)$, and so
$Tb_i=(i+1)b_{i+1}$.  Thus
\[f(a) = \sum_{i\geq 0} \frac{T^i}{i!}b_0 x^i = e^{Ty}b_0,\]
and so we have an isomorphism of $G$\h modules,
\[\Multi^K_G(A;B) \isom \Hom(A,B).\]
\end{example}

\begin{example}\label{ex:getting_vertexops}
With $G$ taken to be the classical vertex group, the singular multilinear
maps from $A_1\tens A_2$ to $B$ are the $\G{2}$\h invariant maps from
$A_1\tens A_2$ to $B\ps{x,y}[(x-y)^{-1}]$.  Given any such map, say $f$, the
$\G{2}$\h invariance just says that for any $a \tens b \in A \tens B$,
\begin{eqnarray*}
f(Ta\tens b) &=& \partial_x f(a\tens b)\\
f(a\tens Tb) &=& \partial_y f(a\tens b).
\end{eqnarray*}
If we add the additional requirement that $f$ be $G$\h invariant, then we
have 
\begin{eqnarray*}
Tf(a\tens b) &=& f(Ta\tens b) + f(a\tens Tb)\\
&& \partial_x f(a\tens b) + \partial_y f(a\tens b),
\end{eqnarray*}
so $G$ acts on $B\ps{x,y}[(x-y)^{-1}]$ as $T_B = \partial_x + \partial_y$.
Taking $V = A_1 = A_2 = B$, and setting
\[Y(\cdot, x)\cdot = f(\cdot \tens \cdot)|_{y=0}:V\tens V \longrightarrow
V\ps{x}[x^{-1}],\]  
we have a vertex operator that satisfies the translation covariance axiom of
equation \eqref{e:inv_axiom}.  
\end{example}

\begin{example}\label{ex:getting_singmaps}
Having seen how to compute a vertex operator from a singular map, we now show
how to compute a binary singular map from a vertex operator.  Let
$Y(\cdot,x)\cdot$ be a vertex operator on a complex vector space $V$.  We
claim that the map, 
\[f(\cdot \tens \cdot)= e^{Ty}Y(\cdot,x-y)\cdot,\] 
gives the desired $G$\h invariant binary singular map.  (Here the operator
$e^{Ty}=\sum_{i\geq 0} y^i \D{i}$ provides a linear map from $V$ to
$V\ps{y}$, and we are regarding $Y(\cdot,x-y) \in V\ps{x-y}[(x-y)^{-1}]$ as
an element of $V\ps{x,y}[(x-y)^{-1}]$ under binomial expansion.  We shall see
in example \ref{ex:binary_with_tail} that this inclusion is important for the
general description that will follow.)  For any $a\tens b \in V\tens V$,
clearly we see that $f(a\tens b)|_{y=0} = Y(a,x)b$, and
\begin{eqnarray*}
f(a\tens b)|_{x=0} &=& e^{Ty}Y(a, -y)b \\
&=& Y(b,y)a
\end{eqnarray*}
where the second equality follows by quasisymmetry (see appendix
\ref{a:axiom_facts}).  Notice also that we could have just as easily
reconstructed $f(a\tens b)$ from $e^{Tx}Y(b,y-x)a$.  Thus we see that a
binary singular map incorporates the quasisymmetry of a vertex operator.

This example provides an interesting contrast with the case for holomorphic
vertex algebras in section \ref{ss:extended}.  There we saw that the entire
theory was determined by a vertex operator $Y(\cdot,x)\cdot$ evaluated at
$x=0$.  Because we have singularities here, we can not make such an
evaluation, but we still can use the actions of $G$ to reconstruct a singular
map from the vertex operator.
\end{example}

The next obvious question is what happens when we compose two singular
multilinear maps?  We begin with an example using the classical vertex group:

\begin{example}
Following the discussion of composition in section \ref{ss:extended}, we
shall assume that all maps are $G$\h invariant, except possibly the bottom
map.  With an eye to emulating the behaviour of axiomatic vertex algebras, we
shall also consider pointwise composition of these maps.

Taking the classical vertex group, $G$, the simplest nontrivial example of
composition is the composition of $f \in \Multi^K_G(A_1, A_2;B_1)$ with $h
\in \Multi^K(B_1, B_2; C)$.  Writing these out explicitly, we have:
\begin{align}
	\begin{aligned}
	A_1 \tens A_2 & \xrightarrow{f} B_1\ps{x_1, x_2}[(x_1-x_2)^{-1}]
	\end{aligned}
				&&
	&\begin{aligned}
	B_1 \tens B_2 & \xrightarrow{h} C\ps{z_1, z_2}[(z_1-z_2)^{-1}].
	\end{aligned}
\end{align}
These compose to give a map 
\[A_1 \tens A_2 \tens B_2 \xrightarrow{h \circ f} C\ps{z_1,
z_2}[(z_1-z_2)^{-1}]\ps{x_1, x_2}[(x_1-x_2)^{-1}].\] 
By the $G$\h invariance of $f$, we know that this $\G{2}$\h linear map
satisfies $\partial_{z_1} = \partial_{x_1} + \partial_{x_2}$.  Using these
differential equations, the composite $h \circ f$ can be seen to factor
through $C\ps{X_1,X_3}[(X_1-X_3)^{-1}]\ps{X_1-X_2}[(X_1-X_2)^{-1}]$ under
either of the following changes of variables:
\begin{align}
	X_1 & = x_1 + z_1   &   X_1 & = x_2 + z_1 \label{e:varchange1}\\
	X_2 & = x_2 + z_1   &   X_2 & = x_1 + z_1 \label{e:varchange2}\\
	X_3 & = z_2         &   X_3 & = z_2.      \label{e:varchange3}
\end{align}
\end{example}

We notice immediately that even though the composite factors through a
collection which has only 3 variables, it does not factor through
$\Fun(\G{3},B) \isom B\ps{x_1, x_2, x_3}[(x_1-x_2)^{-1}, (x_1-x_3)^{-1},
(x_2-x_3)^{-1}]$.  In fact it is simple to find examples of such composites
which do not fit inside $\Fun(\G{3},B)$.  Thus our singular multilinear maps
do not compose in the way we might have initially hoped.

Since we haven't yet defined singular multilinear maps for three input
modules, $\Multi^K(A_1, A_2, A_3;B)$, we could think of defining it to be a
space large enough to contain the composite described in the previous
example, and all other possible composites.  But resulting spaces would end
up simply being the spaces of formal distributions which we had initially set
out to avoid.  Instead, we shall define a collection of singular multilinear
maps for every way of composing singular multilinear maps.  In other words,
if we denote the collection $\Multi^K(A_1, A_2;B)$ by the labelled tree
\[\flattwoleafed{A_1}{A_2}{B}{}{}\]
then the composite would be an element of a space associated to the labelled
tree 
\[\threeleafedleft{A_1}{A_2}{B_2}{C}{}{}{}{}\]
and so to every binary tree $p$, we would associate a collection of singular
multilinear maps, $\Multi_{p}^K(A_1, \ldots, A_n;C)$.  How then shall we
define $\Multi_{\epsfig{file=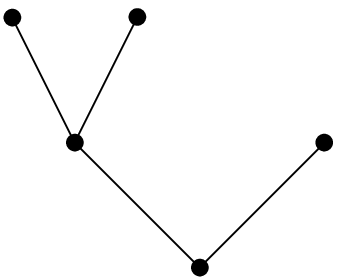,height=3mm}}^K(A_1, A_2,
A_3;C)$?  Following the discussion of section \ref{ss:motivation} we shall
take the following definition:

\begin{dfn}
The collection of multilinear singular maps from $A_1, A_2, A_3$ to $C$
associated to the tree \epsfig{file=Images/left_3.eps,height=4mm} is given by
the following equation:
\[\Multi_{\epsfig{file=Images/left_3.eps,height=3mm}}^K(A_1, A_2,
A_3;C) := \Hom\Bigl(A_1\tens A_2, K \tens \Hom\bigl(\G{2}\tens A_3, K \tens \Hom(\G{2},
C)\bigr)\Bigr)/\sim,\]
where we are quotienting by four actions of $G$.  These will be clear by
regarding this collection of multi maps inside the larger collection:
\[\Hom_{\G{3}}\Bigl(A_1\tens A_2\tens A_3, K \tens \Hom_{G}\bigl(\G{2}, K \tens \Hom(\G{2},
C)\bigr)\Bigr).\]
Then as with the binary singular maps described above, the action on each of
the $A_i$ carries over to an action on its respective occurrence of $G$, and
the additional action of $G$ acts on the one hand on the remaining occurrence
of $G$, and on the other hand it acts diagonally on the outer $\G{2}$ term.
\end{dfn}

\begin{note}
Would it be better to use something like this:
\begin{eqnarray*}
\Multi_{\epsfig{file=Images/left_3.eps,height=3mm}}^K(A_1, A_2,
A_3;C) &:=&  \Hom_{\G{3}}(A_1\tens A_2\tens R, K \tens \Hom_{\G{2}}(\G{2}\tens
A_3, K \tens \Hom(\G{2}, C)))\\
&=&  \Hom_{\G{2}}(A_1\tens A_2, \Fun_{\G{2}}(\G{2}\tens A_3, \Fun(\G{2}, C))).
\end{eqnarray*}
\end{note}

\begin{example}
Taking $G$ to be the classical vertex group, this just says that the
collection of singular multilinear maps,
$\Multi_{\epsfig{file=Images/left_3.eps,height=3mm}}^K(A_1, A_2, A_3;C)$ is
the subcollection of 
\[\Hom\Bigl(A_1\tens A_2, \Hom\bigl(A_3, C\ps{x_1,
x_2}[(x_1-x_2)^{-1}]\bigr) \ps{y_1, y_2}[(y_1-y_2)^{-1}]\Bigr),\] 
satisfying the equations 
\begin{eqnarray*}
T_{A_1} &=& \partial_{y_1}, \\
T_{A_2} &=& \partial_{y_2},\\
T_{A_3} &=& \partial_{x_2}\text{, and}\\   
\partial_{y_1}+\partial_{y_2} &=& \partial_{x_1},
\end{eqnarray*}
where $T_{A_i}$ denotes the action of $G$ on the module, $A_i$.  We can
represent this subcollection pictorially by the labelled tree
\[\threeleafedleft{A_1}{A_2}{A_3}{C}{y_1}{y_2}{x_1}{x_2}\]
or consider it to be the appropriate subcollection of
\[\Hom\Bigl(A_1\tens A_2\tens A_3, C\ps{x_1,x_2}[(x_1-x_2)^{-1}] \ps{y_1,
y_2}[(y_1-y_2)^{-1}]\Bigr),\]  
where the singularity $(y_1-y_2)^{-1}$ is ``independent'' of $A_3$.
\end{example}

We can extend this definition naturally to include all binary trees that
have exactly one splitting at each level.  A difficulty arises
when we consider two compositions at a single level.  For example, the
collection of singular maps associated to the binary labelled tree:
\[\fourleafdoubleB{A_1}{A_2}{A_3}{A_4}{C}\]
We could arrive at a singular multilinear map of type
\epsfig{file=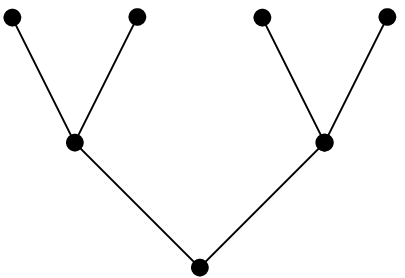,height=3mm} via composition in three ways: for
$G$\h modules $B_1$ and $B_2$, we could compose three binary singular maps:
\begin{equation} \label{e:sing1}
\Multibin(A_1,A_2;B_1) \tens \Multibin(A_3,A_4;B_2) \tens \Multibin(B_1,B_2;C) 
\longrightarrow 
\Multi^K_{\epsfig{file=Images/4ht2.eps,height=3mm}}(A_1, A_2, A_3, A_4;C),
\end{equation}
or we could compose an appropriate binary singular map with a singular map of
type \epsfig{file=Images/left_3.eps,height=3mm} or
\epsfig{file=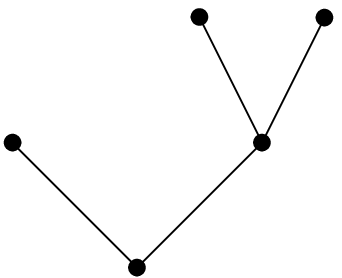,height=3mm} as:
\begin{eqnarray}
\Multi_{\epsfig{file=Images/left_3.eps,height=3mm}}^K(A_1, A_2, B_2;C)
\tens \Multibin(A_3,A_3;B_2) &\longrightarrow & 
\Multi^K_{\epsfig{file=Images/4ht2.eps,height=3mm}}(A_1, A_2, A_3, A_4;C) \label{e:sing2}\\
\Multibin(A_1,A_2;B_1) \tens
\Multi_{\epsfig{file=Images/right_3.eps,height=3mm}}^K(B_1, A_3, A_4;C)
&\longrightarrow &  
\Multi^K_{\epsfig{file=Images/4ht2.eps,height=3mm}}(A_1, A_2, A_3, A_4;C).\label{e:sing3}
\end{eqnarray}
We would like to end up with a space which consists exactly of maps that could
arise as composites.  The important feature of this definition is the
restriction of the dependence of each singularity to the relevant controlling
modules.  This suggests the following definition:

\begin{dfn}\label{d:4leaf_height2}
The collection of singular multilinear maps of type
\epsfig{file=Images/4ht2.eps,height=3mm} from $A_1, A_2, A_3, A_4$ to $C$ is
defined to be the pullback of the singular multilinear maps,
\begin{eqnarray}
\Hom_{G,\G{2}}\Bigl(A_1\tens A_2,K\tens \Hom_{G,\G{2}}\bigl(A_3\tens A_4,
K\tens \Hom_{G,\G{2}}(\G{4}, K\tens\Hom(\G{2},C))\bigr)\Bigr)
\label{e:4multi_1}\\ 
\Hom_{G,\G{2}}\Bigl(A_3\tens A_4,K\tens \Hom_{G,\G{2}}\bigl(A_1\tens A_2,
K\tens \Hom_{G,\G{2}}(\G{4}, K\tens\Hom(\G{2},C))\bigr)\Bigr)
\label{e:4multi_2}
\end{eqnarray}
over 
\begin{equation}
\Hom_{\G{4}}\Bigl(A_1\tens A_2\tens A_3\tens A_4, K\tens K\tens \Hom_{\G{2}}\bigl(\G{4},
K\tens\Hom(\G{2},C)\bigr)\Bigr)
\end{equation}
where the appropriate $G$\h invariance is taken into account.
\end{dfn}

\begin{note}
Because invariance just appears as an equaliser, it commutes with pullbacks.
\end{note}

\begin{thm}\label{t:4compose}
The composites of singular maps in equations \eqref{e:sing1},
\eqref{e:sing2}, and \eqref{e:sing3} compose as desired. 
\end{thm}

Before we prove this theorem, we need the following lemma.

\begin{lemma}
In any symmetric monoidal category, $\cC$, the following diagram commutes:
\[\xymatrix{ \Hom(A_1, B_1\tens C_1)\tens\Hom(A_2,B_2\tens C_2)\ar[r]\ar[d]& \Hom(A_1,
C_1 \tens \Hom(A_2,B_1\tens B_2\tens C_2))\ar[d]  \\
\Hom(A_2, C_2 \tens \Hom(A_1,B_1\tens C_1\tens B_2))\ar[r]&
\Hom(A_1\tens A_2, B_1\tens C_1\tens B_2\tens C_2)}\]
\end{lemma}

\begin{proof}
The proof follows immediately from the fact that the evaluation of 
\[\Hom(A_1, B_1\tens C_1)\tens\Hom(A_2,B_2\tens C_2)\] 
on $A_1\tens A_2$ gives the same result when carried out by either first
evaluating $A_1$, or by first evaluating $A_2$ or by evaluating both
together.
\end{proof}

\begin{proof}[Proof of theorem \ref{t:4compose}]
We shall only prove that the composition works for the classical vertex
group.  For the general theory, see \cite{cts_valgdetails}.  Since
composition is pointwise, the composite \eqref{e:sing1} factors through both
\eqref{e:sing2} and \eqref{e:sing3}, and so we can focus on those
compositions.  In order to prove that the composites map into the pullback,
we shall first prove that they map into each of the pullback objects (the
singular maps given in \eqref{e:4multi_1} and \eqref{e:4multi_2}).

Taking the composite in equation \eqref{e:sing2}, we see that it can be
considered a multi map as in equation \eqref{e:4multi_1} through the following
natural map:
\begin{equation*}
\begin{split}
\Multi&_{\epsfig{file=Images/left_3.eps,height=3mm}}^K(A_1, A_2, B_2;C)
\tens \Multibin(A_3,A_3;B_2)=\\
\Hom&\Bigl(A_1\tens A_2, K \tens \Hom\bigl(\G{2}\tens B_2, K \tens
\Hom(\G{2}, C)\bigr)\Bigr) \tens \Hom\Bigl(A_3\tens A_4, K \tens
\Hom(\G{2},B_2)\Bigr)\\
\longrightarrow& \Hom\biggl(A_1\tens A_2, K \tens \Hom\bigl(\G{2}\tens B_2, K \tens
\Hom(\G{2}, C)\bigr) \tens \Hom\Bigl(A_3\tens A_4, K \tens
\Hom(\G{2},B_2)\Bigr)\biggr) \\
\longrightarrow& \Hom\biggl(A_1\tens A_2, K \tens \Hom\Bigl(A_3\tens A_4, K \tens
\Hom(\G{2},\Hom\bigl(\G{2}, K \tens
\Hom(\G{2}, C)\bigr) )\Bigr)\biggr) \\
\isom & \Hom\biggl(A_1\tens A_2, K \tens \Hom\Bigl(A_3\tens A_4, K \tens
\Hom(\G{4}, K \tens \Hom(\G{2}, C))\Bigr)\biggr).
\end{split}
\end{equation*}
Similarly we can can see that it gives a multi map as in equation
\eqref{e:4multi_2} through evaluating $A_3 \tens A_4$ first.  From the lemma,
we know that these two ways of evaluating are equal, so the composite must
factor through the pullback.
\end{proof}

In the next section we give a general definition for singular multilinear
maps parameterised by any binary tree.  But before we do so, we finish this
section with a description of
$\Multi^K_{\epsfig{file=Images/4ht2.eps,height=3mm}}(A_1, A_2, A_3, A_4;C)$
for the classical vertex group.

\begin{example}
Letting $G$ be the classical vertex group, the collection of singular
multilinear maps of type \epsfig{file=Images/4ht2.eps,height=3mm} from $A_1,
A_2, A_3, A_4$ to $C$ is the submodule of 
\[\Hom\Bigl(A_1\tens A_2\tens A_3\tens
A_4,C\ps{x_1,x_2}[(x_1-x_2)^{-1}] \ps{y_1,y_2, z_1,
z_2}[(y_1-y_2)^{-1},(z_1-z_2)^{-1}]\Bigr),\]
where the singularity $(y_1-y_2)^{-1}$ is independent of $A_3\tens A_4$,
the singularity $(z_1-z_2)^{-1}$ is independent of $A_1 \tens A_2$, and the
maps satisfy the following equations:
\begin{align*}
T_{A_1}&=\partial_{y_1} & T_{A_2}&=\partial_{y_2}\\ 
T_{A_3}&=\partial_{z_1} & T_{A_4}&=\partial_{z_2}\\ 
\partial_{y_1}+\partial_{y_2} &= \partial_{x_1} &
\partial_{z_1}+\partial_{z_2} &= \partial_{x_2} 
\end{align*}
We represent this collection pictorially by the following labelled tree:
\[\fourleafdoubleC{A}{C}{y}{z}{x}\]
The $G$\h invariant such maps are the further subcollection satisfying
$T_C=\partial_{x_1}+\partial_{x_2}$.  
\end{example}

\subsection{Multi Maps Parameterised by Binary Trees}\label{ss:multimaps_binary}
In this section we are concerned with giving the definition of a multilinear
singular map associated to an arbritrary binary tree.  Given any binary tree
$p$, we may consider its collection of internal vertices.  For our purposes,
we shall assume that these include the root, but they do not include the
leaves.  Considering them as a set, this set inherits a partial order from
the tree, where the root is the least element.  We know that any partial
order can be extended to at least one total ordering, possibly many.

We have already been working with trees whose leaves and root are labelled by
$G$\h modules.  It will be useful for the explanations to follow to assume
that every internal node is also labelled.  For any internal node, $q$,
connected to $n$ \defn{incoming} nodes (i.e., non-empty nodes whose height is
equal to the height of $q$ plus one and connected to $q$ by a single edge),
we shall label $q$ by $\G{n}$.  We can also associate to $q$ the tensor
product of the labels of the incoming nodes, and denote it $X_q$.  Thus the
following labelled tree has two internal nodes,
\[\threeleafedleft{A_1}{A_2}{A_3}{C}{}{}{}{}\]
and we have $X_{\text{root}}=\G{2}\tens A_3$, and
$X_{\text{internal}}=A_1\tens A_2$.

We say that our tree is \defn{augmented} when we have added an additional
vertex and edge to the root of the tree, such that the new vertex inherits
the former root label.  The former root is labelled $\G{2}$ as expected.  We
denote the new vertex $\bottom$, and we automatically have $X_\bottom =
\G{2}$.

\begin{dfn}\label{d:bin_sing_maps}
Let $p$ be a binary $n$\h leafed tree, and let $t$ denote a total ordering,
$\bottom < root < p_1 < \cdots < p_l$, of the internal vertices of augmented
$p$, compatible with the the partial ordering inherited from the tree
structure of $p$.  We define an operator on $\G{2}$\h modules:
\[\Sing_{p_i} = \Hom_{G,\G{2}}(X_{p_i}, K\tens\cdot).\]
Iterating this operator we have
\[ \mathrm{Ord}_t(A_1, \ldots, A_n, C) = \Sing_{p_l}\cdots \Sing_{p_1}
\Sing_{\text{root}} \Hom(X_\bottom,C).\] 
Then $\Multi^K_p(A_1, \ldots, A_n, C)$ is defined to be the wide pullback of each
$\mathrm{Ord}_t$ for all possible total orderings, $t$, of of the internal vertices
of augmented $p$, over 
\[\Hom_{\G{n}}(A_1\tens\cdots\tens A_n, \Fun_{p}(\G{n},C))\]
where $\Fun_{p}(\G{n},C))$ is the collection of singular functions as defined
in \cite{bor}.  See appendix \ref{ss:bor_sing_maps} for more details.

We have \defn{$G$\h invariant multilinear singular maps}, $\Multi^K_{G,p}$
exactly when $\Hom(X_\bottom,C)$ is $G$\h invariant.
\end{dfn}

\begin{example}
If $p$ is a tree with only one binary tree at each level, then there is only
one total ordering, $t$, of internal vertices of the tree, and so 
\[\Multi^K_p(A_1,\ldots, A_n, C) = \mathrm{Ord}_t(A_1, \ldots, A_n, C).\]
\end{example}

\begin{example}
When $p = \epsfig{file=Images/4ht2.eps,height=3mm}$, there are exactly two
total orderings of internal vertices of this tree, and the corresponding
$\mathrm{Ord}_t(A_1, \ldots, A_4, C)$ functions are given by
\begin{eqnarray}
\mathrm{Ord}_{t_1} &=& \Hom\Bigl(A_1\tens A_2,K\tens
\Hom\bigl(A_3\tens A_4,K\tens \Hom(\G{4},
K\tens\Hom(\G{2},C))\bigr)\Bigr)\\  
\mathrm{Ord}_{t_2} &=& \Hom\Bigl(A_3\tens A_4,K\tens
\Hom\bigl(A_1\tens A_2, K\tens \Hom(\G{4},
K\tens\Hom(\G{2},C))\bigr)\Bigr),
\end{eqnarray}
where the $\Hom$ functions are both $G$\h invariant and $\G{2}$\h invariant.
Since we are pulling back over 
\begin{multline}
\Hom_{\G{4}}(A_1\tens\cdots\tens A_4, \Fun_{p}(\G{4},C)) = \\
\Hom_{\G{4}}(A_1\tens\cdots\tens A_4, K\tens K\tens \Hom_{\G{2}}(\G{4},
K\tens\Hom(\G{2}, C))), 
\end{multline}
this definition reduces to definition \ref{d:4leaf_height2}.
\end{example}

\section{Vacuum and Locality - Extending from Binary Trees}
Now that we have defined spaces of (generalised) vertex operators and their
composites, we shall return our attention to the remaining axioms for a
vertex algebra.  In particular, we will be interested in incorporating the
vacuum axioms of equations \eqref{e:vac_triv_axiom}, \eqref{e:vac_id_axiom},
and \eqref{e:vac_ps_axiom} into our description, and making sense of the
locality axiom of equation \eqref{e:locality}.  But first, we return to the
case of singular maps between $G$\h modules $A$ and $B$, and parameterise
them by non-branching trees.

\subsection{Non-branching Trees}\label{ss:non_branch_trees}
When we first considered singular multilinear maps in section
\ref{ss:sing_mult_maps}, we defined $\Multi^K(A;B)$ to be
$\Hom_G(A,\Hom(G,B))$ (see definition \ref{d:sing_1and2}).  Now that we have
singular multilinear maps parameterised by trees, we shall think of this
collection of maps as being associated to the flat tree with one leaf,
\epsfig{file=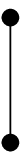,height=3mm}.  (For a review of definitions
related to trees, see appendix \ref{a:cat_of_trees}  Thus we have 
\[\Multi^K_{\ \epsfig{file=Images/1_flat.eps,height=3mm}}(A;B) = \Multi^K(A;B)
= \Hom_G(A,\Hom(G,B)).\]

In fact, we can extend the construction of $\Multi^K_{p}$, described in
definition \ref{d:bin_sing_maps} for binary trees, to trees with
non-branching subtrees.  We do this by simply extending the definition of
$\Sing_{q}$ to those internal nodes $q$ with only one incoming edge. 

\begin{dfn} \label{d:1_sing_maps}
If the tree $p$ in definition \ref{d:bin_sing_maps} is allowed to also have
non-branching subtrees, then $\Multi^K_p(A_1, \ldots, A_n;C)$ is defined
exactly as in that definition except that when an internal vertex, $p_i$ has
only one incoming edge, we define an operator to act on $G$\h modules,
\[\Sing_{p_i} = \Hom_G(X_{p_i},\cdot),\]
where $X_{p_i}$ is the label of the incoming node as in section
\ref{ss:multimaps_binary}.  
\end{dfn}

\begin{lemma}
The definition of $\Multi^K_{\ \epsfig{file=Images/1_flat.eps,height=3mm}}(A;B)$
coincides with definition \ref{d:1_sing_maps} applied to the tree
\epsfig{file=Images/1_flat.eps,height=3mm}. 
\end{lemma}

\begin{proof}
Repeating the construction of definition \ref{d:1_sing_maps} for the tree
\epsfig{file=Images/1_flat.eps,height=3mm}, we have
\[\Multi^K_{\ \epsfig{file=Images/1_flat.eps,height=3mm}}(A;B) =
\Hom_G(A, \Hom(G, B)) \isom \Hom(A,B) \] 
as defined above.  We also see that the $G$\h
invariant maps are the collection,
\[\Multi^K_{G,\epsfig{file=Images/1_flat.eps,height=3mm}}(A;B) 
= \Hom_G(A,\Hom_G(G,B)) \isom \Hom_G(A,B).\]\end{proof}

\begin{example}
If $p$ is the non-branching tree with $n$ internal nodes (including the root)
then 
\[\Multi^K_p(A;B) =
\Hom_G(A,\bigl(\Sing_{\ \epsfig{file=Images/1_flat.eps,height=3mm}}\bigr)^{(n-1)}
\Hom(G,B)).\]
Taking advantage of the $G$\h invariance of 
$\Sing_{\ \epsfig{file=Images/1_flat.eps,height=3mm}}$, we see that this
collection is isomorphic to $\Hom_G(A,\Hom(G,B))$ in the usual way.
\end{example}

\begin{example}
If we take $G$ to be the classical vertex group, then 
$\Multi^K_{\ \epsfig{file=Images/1_flat.eps,height=2mm}}(A;B) \isom
\Hom_{G}(A, B\ps{x}).$ (c.f., example \ref{ex:sing_1_hom}.)  

For the non-branching tree with $n$ internal nodes, $p$, we have  
$\Multi^K_p(A;B) \isom \Hom(A, B\ps{x_1, \ldots, x_n}),$ where the $G$\h
invariance of $\Sing_{\ \epsfig{file=Images/1_flat.eps,height=3mm}}$ means
that any such map is an element of the subcollection, $\Hom(A, B\ps{x_1+
\cdots + x_n})$. 
\end{example}

\begin{example}\label{ex:binary_with_tail}
Consider the singular maps associated to the tree,
\epsfig{file=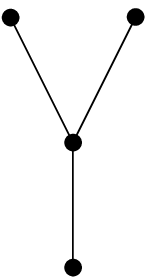,height=3mm}.  From the definition we have
\begin{eqnarray*}
\Multi^K_{\epsfig{file=Images/2_tail.eps,height=3mm}}(A_1,A_2;B) 
&=& \Hom_{G,\G{2}}\Bigl(A_1\tens A_2, K\tens\Hom_G\bigl(\G{2},\Hom(G,B)\bigr)\Bigr)\\
&\isom& \Hom_{G,\G{2}}\Bigl(A_1\tens A_2, K\tens\Hom\bigl(\G{2},B\bigr)\Bigr).
\end{eqnarray*}
When $G$ is the classical vertex group this says that there is a bijection
between the collection of multilinear singular maps associated to the
following trees:
\begin{equation}\label{e:binary_with_tail}
\twoleafedtail{A_1}{A_2}{B}{x}{y}{z} \isom \flattwoleafed{A_1}{A_2}{B}{x+z}{y+z}
\end{equation}
This isomorphism follows immediately from the $G$\h invariance at the
internal node, where it provides the relation $\partial_x +\partial_y =
\partial_z$.  The map between these spaces is just given by binomial
expansion as in example \ref{ex:getting_singmaps}, and the singularity
remains unaffected.  When these maps are fully $G$\h invariant, the single
edge attached to the bottom of the tree corresponds to operating on each
singular map with
\[e^{Tz} = \sum_{i\geq 0} \D{i}z^i:B \longrightarrow B\ps{z}.\]

\end{example}

We finish this section by noting in passing that composing any tree with the
tree consisting of a single node, $\bullet$, leaves the tree unchanged.  So
we would like to define $\Multi^K_\bullet$ so that it composes with a
singular map of type $p$ (for some tree, $p$) to give a singular map of type
$p$.  This suggests the following definition:

\begin{dfn}\label{d:id_sing_maps}
For any $G$\h module, $A$, the singular multilinear maps associated to the
tree $\bullet$ are just the endomorphisms of $A$.
\[\Multi^K_\bullet(A;A) = \Hom(A,A).\]
\end{dfn}

\begin{lemma}
The definition of $\Multi^K_\bullet(A;A)$ coincides with definition
\ref{d:1_sing_maps} applied to the tree $\bullet$.
\end{lemma}

\begin{proof}
Definition \ref{d:1_sing_maps} applied to $\bullet$ gives
$\Multi^K_\bullet(A;A) = \Hom(A,A))$ as desired.
\end{proof}

\subsection{Vacuum - Trees with no Leaves}\label{ss:vacuum}
Now that we have a description of unary singular maps, it is natural to
consider the nullary type multi maps.  By way of motivation, recall that the
vacuum was defined to be a distinguished vector denoted $\vac \in V$ such
that for any $a \in V$, the vacuum satisfies:
\begin{eqnarray*}
T\vac &=& 0 \\
Y(\vac, x)a &=& a \\
Y(a,x)\vac|_{x=0} &=& a.
\end{eqnarray*}
We saw in section \ref{ss:extended} that for vertex algebras without
singularities, the vacuum acted as a unit object.  Later we shall show that
even with singularities, the vacuum will continue to act as a unit.  But in
this section we are concerned primarily with showing how the collection of all
possible vacuum vectors arises naturally when considering singular maps
parameterised by trees.

They arise naturally as the multilinear singular map associated to the empty
tree (which we denote $\circ$).  Applying definition \ref{d:1_sing_maps}, to
the case of the empty tree, we construct the collection of multilinear
singular maps by first augmenting the empty tree, giving
\epsfig{file=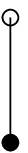,height=3mm}.  As with the tree $\bullet$, the
only internal node is $\bottom$.  Thus we have $X_\bottom = R$ since it is
not connected to any other leaf or internal node (the empty node is not
counted as a leaf), and so we have:

\begin{dfn}\label{d:0_sing_maps}
The multilinear singular maps parameterised by the empty tree are given by:
\[ \Multi^K_{\circ}(R;V) = \Hom(R,V) \isom V.\]
\end{dfn}
The $G$\h invariant elements of this collection are the vectors in $V$ with
trivial $G$\h action, so we see that a $G$\h invariant singular multilinear
map of type $\circ$ is just a vector $v \in V$ satisfying the first vacuum
axiom.

What happens when the vacuum vector is used as an input for a binary singular
multi map?  We begin with an example for the classical vertex group.

\begin{example}\label{ex:compose_with_vac}
Let $G$ be the classical vertex group and $f$ be a binary singular map in 
\[\Multibin(A_1,A_2;B) \isom \Hom_{\G{2}}(A_1\tens A_2, B\ps{x,y}[(x-y)^{-1}]).\] 
We know that for a vector $v \in A_1$ arising as above, we have $Tv=0$, so
given any $a \in A_2$ we have,
\[\partial_x f(v\tens a)= f(Tv\tens a)=0.\]
So $f(v\tens a) \in B\ps{y}$.  We see further what if we assume $f$ to be
$G$\h invariant, we have that $Tf(v\tens a) = \partial_y f(v\tens a)$,
and so writing $f(v\tens a) = \sum_{i\geq 0} b_i y^i$, we see that
$Tb_i=(i+1)b_{i+1}$, and so we have
\begin{equation}\label{e:exp_T}  
f(v\tens a) = \sum_{i\geq 0} \frac{T^i}{i!}b_0 y^i = e^{Ty}b_0.
\end{equation}
But we saw in the previous section that
$\Multi^K_{\ \epsfig{file=Images/1_flat.eps,height=2mm}}(A_2;B) \isom
\Hom_{G}(A_2, B\ps{y}),$ and the $G$\h invariant maps in this collection were
of the form of equation \eqref{e:exp_T}.  Thus we see that we have a
composition map
\begin{equation}\label{e:null_comp}
\Multibin(A_1,A_2;B) \tens \Multi^K_{G,\circ}(R;A_1) \longrightarrow
\Multi^K_{\ \epsfig{file=Images/1_flat.eps,height=2mm}}(A_2;B),
\end{equation}
and $G$\h invariant binary singular maps compose to give $G$\h invariant maps
$\Multi^K_{G,\epsfig{file=Images/1_flat.eps,height=2mm}}(A_2;B)$.  Notice
that if $A_1 = A_2 = B = V$ for some complex vector space $V$, then equation
\eqref{e:exp_T} gives the remaining vacuum axioms when $b_0 = a$.  We shall
see later that this will hold for a suitable algebra.
\end{example}

This is a very satisfactory notion for composition because we have taken a
binary tree and composed it with the empty node to give a tree with only one
leaf.  In fact, this holds equally well for an arbritrary vertex group, with
equation \eqref{e:null_comp} holding for any $G$\h modules $A_1, A_2$ and
$B$.  For the general theory see \cite{cts_valgdetails}.

\subsection{Locality}\label{ss:locality}
We have attempted to be quite clear throughout this paper as to the
ramifications of our description for the classical vertex group.  We would
now like to apply the locality axiom to the classical vertex group story so
far, and interpret it in the abstract presentation.  Recall that the locality
axiom says that for any $a, b \in V$, there exists some $N \gg 0$ such that
the following holds:
\[(x-y)^N[Y(a,x), Y(b, y)] = 0.\]
In other words, the following two maps are equal:
\begin{eqnarray*}
(x-y)^NY(a,x)Y(b, y)\cdot:V &\longrightarrow&
V\ps{x}[x^{-1}]\ps{y}[y^{-1}]\label{e:loc_xy}\\ 
(x-y)^NY(b,y)Y(a, x)\cdot:V &\longrightarrow&
V\ps{y}[y^{-1}]\ps{x}[x^{-1}]\label{e:loc_yx}. 
\end{eqnarray*}
In light of our previous discussion showing how $Y(a,x)Y(b, y)\cdot$ and
$Y(b,y)Y(a, x)\cdot$ are actually elements of different spaces of singular
maps, this presents us with the question of how to interpret such an equality
of maps.

The usual way to interpret them is by considering them inside the space of
formal power series $V\ps{x, x^{-1},y, y^{-1}}$.  But since both of the
spaces of power series in equations \eqref{e:loc_yx} and \eqref{e:loc_xy} are
properly contained in the larger space of formal power series, they can only
be equal if they map to the intersection $V\ps{y}[y^{-1}]\ps{x}[x^{-1}] \cap
V\ps{x}[x^{-1}]\ps{y}[y^{-1}]$.

\begin{lemma}
$V\ps{y}[y^{-1}]\ps{x}[x^{-1}] \cap V\ps{x}[x^{-1}]\ps{y}[y^{-1}] 
= V\ps{x,y}[x^{-1}, y^{-1}].$
\end{lemma}

\begin{proof}
It is clear that the right hand side is contained in the intersection, so we
need only to prove that an arbitrary element of the intersection is contained
in the power series on the right.  But the only difference between the power
series $V\ps{x}[x^{-1}]\ps{y}[y^{-1}]$, and $V\ps{x,y}[x^{-1}, y^{-1}]$, is
that in the former, polynomials in the variable $x^{-1}$ can exist as
coefficients of the power series in the variable $y$.  But this can not
occur in $V\ps{y}[y^{-1}]\ps{x}[x^{-1}]$, so the intersection is as given.
\end{proof}

So we see that for some $N \gg 0$, $(x-y)^NY(a,x)Y(b, y)\cdot$ and
$(x-y)^NY(b,y)Y(a, x)\cdot$ are equal as maps from $V$ to $V\ps{x,y}[x^{-1},
y^{-1}]$.  Fixing $N$, we denote this map $g(a\tens b\tens \cdot)$.

Considering $g$ as a map to the larger space, $V\ps{x,y}[x^{-1}, y^{-1},
(x-y)^{-1}],$ of power series with the inverse of $(x-y)$ adjoined.  Then the
map, $(x-y)^{-N}g(a\tens b\tens \cdot)$, is well defined on $V$.  Notice first
that this map is not necessarily equal to $Y(a,x)Y(b, y)\cdot$ or $Y(b,y)Y(a,
x)\cdot$, because $(x-y)^{-N}$ does not exist in the codomain of either of
these maps.  But, if we expand $(x-y)^{-N}$ as a power series in either the
variable $x$ or $y$, we have maps
\begin{equation}\label{e:refine_from_three}
\xymatrix{ & V\ps{x,y}[x^{-1},y^{-1},
(x-y)^{-1}] \ar[ld]_{i_{x,y}} \ar[rd]^{i_{y,x}}&\\ 
V\ps{x}[x^{-1}]\ps{y}[y^{-1}] & &V\ps{y}[y^{-1}]\ps{x}[x^{-1}]}
\end{equation}
which are identities except on $(x-y)^{-1}$, where they are given by
\begin{eqnarray*}
i_{x,y}\bigl((x-y)^{-j-1}\bigr) & = & \sum_{n \in \N} \binom{n+j}{j}
x^{-n-j-1} y^n \\
i_{y,x}\bigl((x-y)^{-j-1}\bigr) & = & (-1)^{j+1}\sum_{n \in \N} \binom{n+j}{j}
x^{n} y^{-n-j-1},
\end{eqnarray*}
then we see that under these maps, the singular map $(x-y)^{-N}g(a\tens
b\tens \cdot)$ canonically maps down to the two composites of vertex
operators:
\begin{eqnarray*}
i_{x,y}\bigl((x-y)^{-N}g(a\tens b\tens \cdot)\bigr)&=& Y(a,x)Y(b, y)\cdot\\
i_{y,x}\bigl((x-y)^{-N}g(a\tens b\tens \cdot)\bigr)&=& Y(b,y)Y(a, x)\cdot.
\end{eqnarray*}
So we have constructed an element of the following collection of singular
maps (with appropriate $G$\h invariance taken):
\begin{equation*}
\Hom\biggl(V\tens V, \Hom\Bigl(V, V\ps{x,y}[x^{-1},
y^{-1}]\Bigr)[(x-y)^{-1}]\biggr),
\end{equation*}
where we have taken special care to emphasise that the singularity at $(x-y)$
depends only upon the outer two copies of $V$.  We shall see shortly that
this dependence of singularities on inputs is a crucial feature of the more
general theory.  But first we examine the effect of
the vacuum on these composites.

\begin{example}
Consider the action of $(x-y)^{-N}g(a\tens b\tens \cdot)$ on the vacuum,
$\vac$.  From its definition we have
\begin{eqnarray*}
(x-y)^{-N}g(a\tens b\tens \vac) &=& (x-y)^{-N}(x-y)^NY(a,x)Y(b, y)\vac\\
&=& (x-y)^{-N}(x-y)^NY(b,y)Y(a, x)\vac,
\end{eqnarray*}
and we know from the vacuum axioms that $Y(a, x)\vac$ and $Y(b, y)\vac$ have
no singularities, therefore $(x-y)^{-N}g(a\tens b\tens \vac)$ is an element
of $V\ps{x,y}[(x-y)^{-1}]$.  From our previous discussion we know that it
maps down to the composites
\begin{eqnarray*}
i_{x,y}\bigl((x-y)^{-N}g(a\tens b\tens \vac)\bigr)&=& Y(a,x)e^{Ty}b\\
i_{y,x}\bigl((x-y)^{-N}g(a\tens b\tens \vac)\bigr)&=& Y(b,y)e^{Tx}a,
\end{eqnarray*}
which are elements of $V\ps{x}[x^{-1}]\ps{y}$ and $V\ps{y}[y^{-1}]\ps{x}$
respectively.  Recall that in example \ref{ex:getting_singmaps} we showed how
to construct a singular map from a vertex operator.  

\begin{lemma}
The binary singular power series given by regarding $e^{Ty}Y(a,x-y)b$ as an
element of $V\ps{x,y}[(x-y)^{-1}]$ is equal to $(x-y)^{-N}g(a\tens b\tens
\vac)$.
\end{lemma}

\begin{proof}
Each of these is an element of $V\ps{x,y}[(x-y)^{-1}]$ which agree when
$x=0$ (and when $y=0$).  Because each satisfies $\partial_x +\partial_y =
\partial_V$ they are completely determined by their values at $x=0$ (or
$y=0$), and so are equal.
\end{proof}

\end{example}

Just as we looked at the action of $(x-y)^{-N}g(a\tens b\tens \cdot)$ on the
vacuum, we consider its action on an arbritrary $c \in V$.  From the vacuum
axioms, we know that $c = Y(c,z)\vac|_{z=0}$, so in fact we shall consider
\[(x-y)^{-N}g\Bigl(a\tens b\tens Y(c,z)\vac\Bigr)\Bigr|_{z=0} = 
(x-y)^{-N}(x-y)^NY(a,x)Y(b, y)Y(c,z)\vac|_{z=0}.\] 
Using the locality axiom, we know that for some $M\gg 0$ the following
equality holds for any $a \in V$:
\[(y-z)^MY(a,x)Y(b, y)Y(c,z)\vac = (y-z)^MY(a,x)Y(c,z)Y(b,y)\vac.\]
Setting $z=0$, we see that as an element of
$V\ps{x,y}[(x-y)^{-1},x^{-1},y^{-1}]$, the singularity $y^{-1}$ of
$(x-y)^{-N}g(a\tens b\tens c)$ depends only on $b$ and $c$.  Similarly, the
singularity $x^{-1}$ depends only on $a$ and $c$.  

Recalling how when we defined the binary singular maps, we defined maps
$f(a\tens b)$ such that setting $y=0$ gave $Y(a,x)b$ and setting $y=0$ gave
$Y(b,y)a$, we see that a more natural definition for $g$ would be as a map to
$V\ps{x,y,z}[(x-y)^{-1}, (y-z)^{-1},(x-z)^{-1}]$.  For any $a,b,c \in V$ we
may define it (with $N\gg 0$ as before):
\[(x-y)^{-N}(x-y)^{N}\Bigl(e^{Tz}Y(a,x-z)Y(b, y-z)c\Bigr).\]

Thus we have a ternary function which, according to the dependence of
singularities on inputs discussed above, can we regarded as an element of:
\begin{eqnarray}
\Hom\biggl(V\tens V, \Hom\Bigl(V, V\ps{x,y,z}[(x-z)^{-1},
(y-z)^{-1}]\Bigr)[(x-y)^{-1}]\biggr),\label{e:ternary_1} \\
\Hom\biggl(V\tens V, \Hom\Bigl(V, V\ps{x,y,z}[(x-z)^{-1},
(x-y)^{-1}]\Bigr)[(y-z)^{-1}]\biggr),\label{e:ternary_2} \\
\Hom\biggl(V\tens V, \Hom\Bigl(V, V\ps{x,y,z}[(y-z)^{-1},
(x-y)^{-1}]\Bigr)[(x-z)^{-1}]\biggr).\label{e:ternary_3}
\end{eqnarray}
All three of these collections of singular maps are properly contained in the
larger collection,
\begin{equation}
\Hom\biggl(V\tens V \tens V, V\ps{x,y, z}[(x-z)^{-1}, (y-z)^{-1},
(x-y)^{-1}]\biggr),\label{e:ternary_4}
\end{equation}
where the dependence of each singularity on input is ignored.  So we are led
to define a new collection of multi maps which we associate to the flat three
leafed tree, \epsfig{file=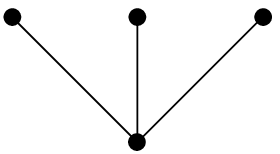,height=3mm}, as the pullback of
the spaces in equations \eqref{e:ternary_1}, \eqref{e:ternary_2}, and
\eqref{e:ternary_3} over the space in equation \eqref{e:ternary_4}.  We shall
denote it $\Multi^K_{\epsfig{file=Images/flat_3.eps,height=3mm}}(V,V,V;V)$.
In fact, this collection formalises the idea of the \defn{operator product
expansion} (see \cite[section 4.6]{kac}).

\begin{note}
We have tacitly been using the isomorphism of trees of equation
\ref{e:binary_with_tail}.
\end{note}

\subsection{Flat Trees and Singular Maps}
In the previous section we used our knowledge of axiomatic vertex algebras to
define a collection of multi maps associated to a three leaved flat tree.
Following that example, we now provide a more general description of singular
maps for an arbritrary vertex group associated to any $n$\h leafed flat tree.
This definition arises naturally from the previous discussion, and
incorporates the definitions provided for trees with $0,1$ and $2$ leaves.

\begin{dfn}\label{d:n_flat_sing_maps}
If $G$ is a vertex group and $A_1, \ldots, A_n, B$ are $G$\h modules, then
the collection of singular maps associated to the flat tree with $n$ leaves,
\epsfig{file=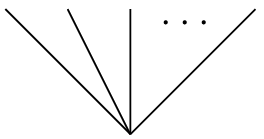,height=3mm}, is denoted
$\Multi^K_{\epsfig{file=Images/generic.eps,height=3mm}}(A_1, \ldots, A_n;
B)$, and is defined to be the pullback of the following singular maps
for $\sigma \in A_n$ (the alternating group):
\begin{multline}\label{e:nary_general_pbk}
\Sing_\sigma(A_1, \ldots, A_n; \Hom(\G{n},B)) =\\ \Hom\biggl(A_{\sigma(1)}\tens
A_{\sigma(2)}, K \tens \Hom\Bigl(A_{\sigma(3)}, K^{\tens 2}\tens \bigl(\cdots
K^{\otimes n-1} \tens \Hom(\G{n},B)\cdots\bigr)\Bigr)\biggr),
\end{multline}
over the collection,
\begin{equation}\label{e:nary_general}
\Hom_{\G{n}}\biggl(A_1\tens \cdots \tens A_n, \Fun(\G{n}, B)\biggr) = 
\Hom_{\G{n}}\biggl(A_1\tens \cdots \tens A_n, K^{\tens \binom{n}{2}} \tens
\Hom(\G{n},B)\biggr),
\end{equation}
where equation \eqref{e:nary_general_pbk} is taken to be invariant under the
action of $\G{n}$ inferred from its action on equation
\eqref{e:nary_general}.  The singularities of equation
\eqref{e:nary_general_pbk} are the singularities are tensored over $f_{i,j}$
as in definition \ref{d:fun}, of $\Fun(\G{n}, B)$.  We have \defn{$G$\h invariant
multilinear singular maps} exactly when $\Hom(\G{n},B)$ is $G$\h invariant.
\end{dfn}

\begin{example}
When $n=0$, this definition says that $\Multi^K_{\circ}(R;B) = \Hom(R,B)$ as
expected, and when $n=1$, we have that $\Multi^K_{\
\epsfig{file=Images/1_flat.eps,height=3mm}}(A;B) = \Hom_G(A,\Hom(G,B))$ as
expected.
\end{example}

\begin{example}
When $n=2$, we see that
$\Multi^K_{\epsfig{file=Images/2_leaf.eps,height=3mm}}(A_1,A_2;B) =
\Hom_{\G{2}}(A_1 \tens A_2, K \tens \Hom(\G{2},B))$ as before.
\end{example}

\begin{example}\label{d:ternary_multimaps}
If $G$ is a vertex group and $A_1, A_2, A_3, B$ are $G$\h modules, then the
collection of singular maps,
$\Multi^K_{\epsfig{file=Images/flat_3.eps,height=3mm}}(A_1, A_2, A_3; B)$, is
the pullback of the following three collections of singular maps:
\begin{eqnarray}
\Hom\biggl(A_1\tens A_2, K \tens \Hom\Bigl(A_3, K^{\tens 2}\tens 
\Hom(\G{3},B)\Bigr)\biggr),\label{e:ternary_general_1} \\ 
\Hom\biggl(A_2\tens A_3, K \tens \Hom\Bigl(A_1, K^{\tens 2}\tens 
\Hom(\G{3},B)\Bigr)\biggr),\label{e:ternary_general_2} \\ 
\Hom\biggl(A_3\tens A_1, K \tens \Hom\Bigl(A_2, K^{\tens 2}\tens 
\Hom(\G{3},B)\Bigr)\biggr),\label{e:ternary_general_3} 
\end{eqnarray}
over the collection,
\begin{equation}
\Hom_{\G{3}}\biggl(A_1\tens A_2 \tens A_3, \Fun(\G{3}, B)\biggr) = 
\Hom_{\G{3}}\biggl(A_1\tens A_2 \tens A_3, K^{\tens 3} \tens
\Hom(\G{3},B)\biggr),\label{e:ternary_general_4}
\end{equation}
with appropriate $G$\h invariance taken.  We see immediately that the
definition of
$\Multi^K_{\epsfig{file=Images/flat_3.eps,height=3mm}}(V,V,V;V)$ for the
classical vertex group follows immediately.
\end{example}

\subsection{Multi Maps Parameterised by All Trees}
After extending our definition of singular maps to trees with many leaves, we
are ready to define the collection of singular maps associated to an
arbitrary tree.  In this section we give the general definition for an
arbritrary vertex group, and show that singular maps compose as desired for
the classical vertex group.  This definition generalises definitions
\ref{d:bin_sing_maps}, \ref{d:1_sing_maps} and \ref{d:n_flat_sing_maps}.  For
clarity we give it in its complete form here.

Recall from section \ref{ss:multimaps_binary} that a tree is said to be
\defn{augmented} when we have added an additional vertex and edge to the root
of the tree, such that the new vertex inherits the former root label, and the
old vertex is labelled $\G{n}$ where $n$ is the number of incoming nodes.  We
denote the new vertex $\bottom$.  Also recall that to each internal node,
$q$, of a tree, we can associate the tensor product of the labels of the
incoming nodes, and denote it $X_q = X_{q,1} \tens \cdots \tens X_{q,n}$.
Thus we have $X_\bottom = \G{n}$.

\begin{dfn}
If $G$ is a vertex group and $A_1, \ldots, A_n, B$ are $G$\h modules, then
the collection of singular maps associated to an arbritrary tree, $p$, with $n$
leaves is denoted, $\Multi^K_p(A_1, \ldots, A_n; B)$, and is defined as
follows:
\begin{itemize}
\item Let $t$ denote a total ordering, $\bottom < \Root < p_1 < \cdots < p_l$,
of the internal vertices of augmented $p$, compatible with the the partial
ordering inherited from the tree structure of $p$.  For each internal vertex,
$p_i$, denote the number of incoming edges $n_i$.  Let $\sigma \in A_{n_i}$
and define an operator taking $\G{n_i}$\h modules to the suitably $\G{n_i}$\h
invariant maps:
\begin{equation}
\Sing_{\sigma,p_i}\cdot= \Hom\biggl(X_{p_i,\sigma(1)}\tens X_{p_i,\sigma(2)}, K
\tens \Hom\Bigl(X_{p_i,\sigma(3)}, K^{\tens 2}\tens \bigl(\cdots K^{\otimes
n_i-1} \tens (\cdot)\cdots\bigr)\Bigr)\biggr).
\end{equation}
The tensor products of the singularities, $K$, are taken over appropriate
copies of $f_{i,j}$ as in definition \ref{d:fun}.
\item Let $\Sing_{p_i}$ denote the pullback of $\Sing_{\sigma,p_i}$
for all $\sigma \in A_{n_i}$ over $\Hom\Bigl(X_{p_i},K^{\otimes \binom{n_i}{2}}
\tens \cdot\Bigr)$.  This can be thought of as giving the multilinear
singular maps associated to the flat subtree with root $p_i$. 
\item Iterating this operator for all internal vertices of $p$, we have
\[ \mathrm{Ord}_t(A_1, \ldots, A_n, B) = \Sing_{p_l}\cdots \Sing_{p_1}
\Sing_{\Root} \Hom(X_\bottom,B).\] 
\end{itemize}
Then $\Multi^K_p(A_1, \ldots, A_n; B)$ is defined to be the wide pullback of each
$\mathrm{Ord}_t$ for all possible total orderings, $t$, of of the internal vertices
of augmented $p$, over 
\begin{equation}\label{e:pull_over_me}
\Hom_{\G{n}}\Bigl(A_1\tens\cdots\tens A_n, \Fun_{p}(\G{n},B)\Bigr)
\end{equation}
where $\Fun_{p}(\G{n},B)$ is the collection of singular functions as defined
in \cite{bor}.  See appendix \ref{ss:bor_sing_maps} for more details.

We have \defn{$G$\h invariant multilinear singular maps}, $\Multi^K_{G,p}$,
exactly when $\Hom(X_\bottom,B)$ is $G$\h invariant.
\end{dfn}

\begin{note}
Because limits commute with one another, we could avoid forming the pullback
$\Sing_{p_i}$ by instead defining
\begin{equation}\label{e:ord_plus_perms}
\mathrm{Ord}_{t,\sigma_I}(A_1, \ldots, A_n, B)  = \Sing_{p_l, \sigma_l}\cdots
\Sing_{p_1, \sigma_1}\Sing_{\Root, \sigma_0} \Hom(X_\bottom,B)
\end{equation}
for each $\sigma_i \in A_{n_i}$, and taking the wide pullback for all
possible total orderings, $t$, and permutations $\sigma_i \in A_{n_i}$ for
$1\leq i \leq l$.
\end{note}

\begin{example}
We see automatically that if $p$ is a flat tree, then there is only one total
ordering, and so $\Multi^K_p(A_1, \ldots, A_n; B)$ reduces to the pullback of
$\Sing_{\Root} \Hom(X_\bottom,B)$, over the maps in equation
\eqref{e:pull_over_me}, which is just definition \ref{d:n_flat_sing_maps}.
\end{example}

\begin{thm}
If $G$ is a vertex group over a field and $\Multi^K_p(A_1, \ldots, A_n; B_1)$,
$\Multi^K_q(B_1, \ldots, B_m; C)$ are collections of multilinear singular
maps, then they compose pointwise to give an element of 
\[\Multi^K_{p\circ q}(A_1,\ldots, A_n, B_2 \ldots, B_m; C).\]
\end{thm}

\begin{note}
Keep in mind that we are composing the trees $p$ and $q$ and not the
augmented trees.  We only use augmented trees for the purpose of describing
their associated singular multi maps.
\end{note}

\begin{proof}
Let $f \in \Multi^K_p$ and $g \in \Multi^K_q$.  We shall show that $f$ and
$g$ compose to give an element of $\Multi^K_{p\circ q}$.  To do so, we select
a pullback object $\mathrm{Ord}_{t,\sigma_I}(A_1, \ldots, A_n, B_2 \ldots,
B_m; C)$ as in equation \eqref{e:ord_plus_perms}.  The total ordering of the
internal vertices of the tree $p \circ q$ provides a total ordering of both
$p$ and $q$ as $\Root_p < p_1 < \cdots < p_k$ and $\Root_q < q_1 < \cdots <
q_l$.  Such a pullback object also consists of a choice of permutation of
labels for each $X_{p_i}, X_{q_j}$, and so we have uniquely determined
objects,
\begin{eqnarray*}
\mathrm{Ord}_{t_p}(A_1, \ldots, A_n, B_1) &=& \Sing_{\sigma_k, p_k}\cdots
\Sing_{\sigma_1, p_1} \Sing_{\sigma_0, \Root_p} \Hom(X_{\bottom_p},B_1)\\
\mathrm{Ord}_{t_q}(B_1, \ldots, B_m, C) &=& \Sing_{\delta_l, q_l}\cdots
\Sing_{\delta_1, q_1} \Sing_{\delta_0, \Root_q} \Hom(X_{\bottom_q},C).
\end{eqnarray*}
Because $f$ and $g$ are contained in the pullbacks, we can regard them as
elements of $\mathrm{Ord}_{t_p}$ and $\mathrm{Ord}_{t_q}$ respectively.

We know that the internal vertices of the composed tree, $p \circ q$ consist
exactly of the internal vertices of $p$, the internal vertices of $q$, and
the root of $p$.  We shall therefore prove that $f$ and $g$ compose to give
an element of $\mathrm{Ord}_{t,\sigma_I}$ by induction on the number of
internal vertices of $p \circ q$, $k+l+1$.  If the composed tree $p \circ q$
has zero internal vertices, then we are composing with a tree of the form
$\bullet$ or $\circ$, and it is clear from the discussions of sections
\ref{ss:non_branch_trees} and \ref{ss:vacuum} that maps compose as desired.

Now, without loss of generality we may assume that in the total ordering $t$,
$p_k > q_l$.  We need only to prove that there exists a natural map,
\begin{equation}\label{e:composition_proof}\begin{split}
\biggl(\Sing_{\sigma_k, p_k}&\cdots \Sing_{\sigma_1, p_1}
\Sing_{\sigma_0,\Root_p} \Hom(X_{\bottom_p},B_1)\biggr) \tens 
\biggl(\Sing_{\delta_l, q_l}\cdots \Sing_{\delta_1, q_1} \Sing_{\delta_0,
\Root_q} \Hom(X_{\bottom_q},C)\biggr)\\  
&\begin{split}\longrightarrow \Sing_{\sigma_k, p_k}\biggl(\Bigl(\Sing_{\sigma_{k-1},
p_{k-1}}&\cdots \Sing_{\sigma_1, p_1} \Sing_{\sigma_0,\Root_p} 
\Hom(X_{\bottom_p},B_1)\Bigr) \tens \\
&\Bigl(\Sing_{\delta_l, q_l}\cdots
\Sing_{\delta_1, q_1} \Sing_{\delta_0, \Root_q} \Hom(X_{\bottom_q},C)\Bigr)\biggr),
\end{split}\end{split}\end{equation}
because on the right hand side, the operator $\Sing_{\sigma_k, p_k}$ is
acting on a collection of singular maps associated to a pair of trees with
$k+l$ internal vertices, which we know compose by induction.  We know that
$\Bigl(\Sing_{\sigma_k, p_k}\cdots \Sing_{\sigma_1, p_1}
\Sing_{\sigma_0,\Root_p} \Hom(X_{\bottom_p},B_1)\Bigr)$ can be evaluated at
$X_{p_k}$ (see appendix \ref{a:evaluation}), so the left hand side of
equation \eqref{e:composition_proof} can be evaluated at $X_{p_k}$, giving
the inner collection of singular maps on the right hand side.  The transpose of
this evaluation map provides us with the desired map.  

Now that we have seen that the maps $f$ and $g$ compose naturally into each
pullback object, $\mathrm{Ord}_{t}$, our proof will be complete if we show
that given two pullback objects, $\mathrm{Ord}_{t_1,\sigma_I}$ and
$\mathrm{Ord}_{t_2,\sigma_J}$, the composite of $f$ and $g$ maps through each
of them to the same element of
\[\Hom\Bigl(A_1\tens \cdots \tens A_n \tens B_2 \tens \cdots \tens B_m,
\Fun_{p\circ q}(\G{n+m-1},C)\Bigr),\] 
the collection over which we are pulling back.  Denote this collection $Z$.
Since both $f$ and $g$ are elements of a pullback over objects 
\begin{eqnarray*}
X &=& \Hom\Bigl(A_1\tens \cdots \tens A_n,\Fun^G_{p}(\G{n},B_1)\Bigr)\\
Y &=& \Hom\Bigl(B_1\tens \cdots \tens B_m,\Fun_{q}(\G{m},C)\Bigr)
\end{eqnarray*}
respectively, then if there exists an injective maps from $Z$ to the
composite of $X$ and $Y$, then we know that the composite of $f$ and $g$ maps
to the same element of $Z$, and hence is an element of the pullback.  Clearly
$X$ and $Y$ compose to give an element of the collection,
\[W = \Hom\Bigl(A_1\tens \cdots \tens A_n \tens B_2 \tens \cdots \tens B_m,
\Fun_{p}\bigr(\G{n},\Fun_{q}(\G{m},C)\bigr)\Bigr),\] 
and because we are working with a vertex group over a field, there is a
natural injection from $Z$ to $W$, so our proof is complete.
\end{proof}

\begin{remark}
Composition can be checked to be associative.
\end{remark}

\subsection{Refinement and Maps Between Singular Functions}\label{ss:refinement}
Now that we have a definition of multilinear singular maps for every possible
tree, we consider how these collections are related.  Recall that in equation
\eqref{e:refine_from_three} we saw that we could map from the collection of
singular maps associated to the tree,
\epsfig{file=Images/flat_3.eps,height=3mm}, to the collections associated to
the trees \epsfig{file=Images/left_3.eps,height=3mm} and
\epsfig{file=Images/right_3.eps,height=3mm} by expanding the singularity,
$(x-y)^{-1}$, as a power series in either $x$ or $y$.  We shall now
generalise this expansion for an arbritrary vertex group and see that this is
an example of a more general phenomenon where our collections of multilinear
singular maps are related to one another through canonical maps.

\begin{dfn}
Given any ring of singular functions, $K$ (see definition \ref{d:K}),
a \defn{refinement for a singularity} is the map,
\[K \longrightarrow \Hom_G(G,K)\]
which takes any $k \in K$ to the map $f \in \Hom_G(G,K)$ defined by $f(g) = g
\acts k$ for any $g \in G$.  For any other $G$\h module, $A$, and $G$\h
invariant map $\alpha:A \rightarrow G$, we define a refinement for $K$ by
composition:
\[K \longrightarrow \Hom_G(G,K) \xrightarrow{\alpha} \Hom_G(A,K). \]
\end{dfn}

We know that this map is well defined since $K$ is a $G$\h module.  It is
also obviously a $G$\h linear map.  The following examples will show that
this is just a generalisation of the idea of power series expansions.

\begin{example}
We saw in section \ref{ss:classicalvg} that $K = \C \ps{x}[x^{-1}]$.  So a
refinement for this singularity is a map:
\begin{eqnarray*}
\C \ps{x}[x^{-1}] &\longrightarrow &\C \ps{x}[x^{-1}]\ps{y} \\
x^k &\mapsto & \sum_{i \geq 0} \binom{k}{i}x^{k-i}y^i.
\end{eqnarray*}
In other words, refinement for this singularity is a map from $x^k$ to
$(x+y)^k$ expanded as a power series in the variable $y$.  
\end{example}

\begin{example}
If we compose the refinement map with the antipode map, $S:G \rightarrow G$,
we have:
\begin{eqnarray*}
\C \ps{x}[x^{-1}] &\longrightarrow &\C \ps{x}[x^{-1}]\ps{y} \\
x^k &\mapsto & \sum_{i \geq 0} \binom{k}{i}(-1)^i x^{k-i} y^i,
\end{eqnarray*}
which is just $(x-y)^k$ expanded as a power series in the variable $y$.
\end{example}

\begin{example}
If we compose the refinement map with the multiplication map, $\mu:\G{2}
\rightarrow G$, we have:
\begin{eqnarray*}
\C \ps{x}[x^{-1}] &\longrightarrow &\C \ps{x}[x^{-1}]\ps{y,z} \\
x^k &\mapsto & \sum_{p,q \geq 0} \binom{k}{p+q} \binom{p+q}{p} x^{k-p-q} y^p z^q,
\end{eqnarray*}
which is just $(x+y+z)^k$ expanded as a power series in the variables $y$ and
$z$.
\end{example}

In order to apply this notion of refinement to our singular multilinear
maps, we shall associate a refinement for a singularity to
every map between trees in the category of trees.  

\begin{thm}\label{t:refinement_works}
For $G$\h module, $A_1, \ldots, A_n, B$, and any tree $p$ with $n$ leaves,
there exists a canonical map $\Multi^K_p(A_1, \ldots, A_n; B) \rightarrow
\Multi^K_q(A_1, \ldots, A_n; B)$ for all trees $q$ which refine to $p$.  This
map is given by appropriate refinements of singularities.
\end{thm}

Before we prove this theorem, we give some examples for the classical vertex
group.

\begin{example}
We saw in example \ref{ex:binary_with_tail} that for the classical vertex
group, there is an isomorphism
\[\twoleafedtail{A_1}{A_2}{B}{x_1}{x_2}{y} \isom
\flattwoleafed{A_1}{A_2}{B}{z_1}{z_2}\] 
\end{example}

\begin{example}
Again for the classical vertex group, consider the refinement map
\[\flattwoleafed{A_1}{A_2}{B}{z_1}{z_2} \longrightarrow
\twoleavedrighttwo{A_1}{A_2}{B}{x_1}{x_2}{y} \]
At the level of power series, this is a map,
\[\Hom\Bigl(A_1\tens A_2, B\ps{z_1, z_2}[(z_1-z_2)^{-1}]\Bigr)
\longrightarrow \Hom\Bigl(A_2, \Hom\bigl(A_1, B\ps{x_1,
x_2}[(x_1-x_2)^{-1}]\bigr)\ps{y}\Bigr) \]
where $z_1$ is mapped to $x_1$, $z_2$ is mapped to $x_2+y$, and
$(z_1-z_2)^{-1}$ is mapped to $(x_1-x_2-y)^{-1}$ and expanded as power series
in the variable $y$.  Also we see that automatically $\partial_y =
\partial_{x_2}$, which is the same as the $G$\h invariance requirement at the
internal node.

Abstractly we have a refinement map
\begin{eqnarray*}
\Hom\Bigl(A_1\tens A_2, \Hom(\G{2},B)\tens K\Bigr)
&\longrightarrow& \Hom\Bigl(A_2, \Hom\bigl(A_1 \tens G, K \tens \Hom(\G{2},B)\bigr)\Bigr) \\
&&\isom \Hom\Bigl(A_1 \tens A_2, \Hom\bigl(G, K \tens \Hom(\G{2},B)\bigr)\Bigr), 
\end{eqnarray*}
and so this map reduces to a map 
\begin{eqnarray*}
\Hom(\G{2},B)\tens K &\longrightarrow &\Hom\bigl(G, K \tens \Hom(\G{2},B)\bigr)\\
f\tens k  &\mapsto & F
\end{eqnarray*}
where $F(g) = \sum_{(g)} f(\cdot \tens g_{(1)}\cdot) \tens (S(g_{(2)})k)$.
\end{example}

\begin{example}\label{ex:expand_three}
If we again let $G$ be the classical vertex group and consider the refinement map,
\[\flatthreeleafedlabelled{A_1}{A_2}{A_3}{B}{z_1}{z_2}{z_3} \longrightarrow  
\threeleafedleft{A_1}{A_2}{A_3}{B}{x_1}{x_2}{y_1}{y_2}\]
At the level of power series, we have a map of variables,
\[B\ps{z_1,z_2,z_3}[(z_1-z_2)^{-1}, (z_1-z_3)^{-1}, (z_2-z_3)^{-1}]
\longrightarrow B\ps{y_1, y_2}[(y_1-y_2)^{-1}]\ps{x_1, x_2}[(x_1-x_2)^{-1}]\] 
where the following are expanded as power series in the variables $x_1, x_2$.
\begin{align*}
z_1 &\mapsto x_1+y_1 & (z_1-z_2)^{-1} &\mapsto (x_1-x_2)^{-1}\\
z_2 &\mapsto x_2+y_1 & (z_1-z_3)^{-1} &\mapsto (x_1+y_1-y_2)^{-1}\\
z_3 &\mapsto y_2     & (z_2-z_3)^{-1} &\mapsto (x_2+y_1-y_2)^{-1}.
\end{align*}

Abstractly, we need to define a map from
$\Multi^K_{\epsfig{file=Images/flat_3.eps,height=3mm}}$ into
\[\Hom\biggl(A_1 \tens A_2, K \tens \Hom\Bigl(\G{2} \tens A_3, K \tens
\Hom(\G{2}, B)\Bigr)\biggr).\]
Of the possible pullback objects corresponding to the flat 3 leafed tree, we
shall show that there exists a canonical map from 
\[\Hom\biggl(A_1 \tens A_2, K \tens \Hom\Bigl(A_3, K^{\otimes 2} \tens
\Hom(\G{3}, B)\Bigr)\biggr).\]
Since the outer terms are the same we need only define a map 
\begin{eqnarray*}
K^{\otimes 2} \tens \Hom(\G{3}, B) &\longrightarrow &\Hom\Bigl(\G{2}, K \tens
\Hom(\G{2}, B)\Bigr)\\
k_1 \tens k_2 \tens f &\mapsto & F,
\end{eqnarray*}
where $F(g\tens h) = \sum_{(g)}\sum_{(h)} \biggl[\Bigl(S(g_{(1)})k_1\Bigr)
\Bigl(S(h_{(1)})k_2\Bigr)\biggr] \tens \biggl[f\Bigl(g_{(2)}\cdot \tens
h_{(2)}\cdot \tens \cdot \Bigr)\circ(\Delta \tens 1)\biggr].$
\end{example}

\begin{example}
From the discussion of locality in section \ref{ss:locality}, we see that by
setting $z_3$, $x_2$ and $y_2$ equal to zero, the refinement map between
$G$\h invariant singular functions,
\[\flatthreeleafedlabelled{V}{V}{V}{}{z_1}{z_2}{0} \longrightarrow  
\threeleafedright{V}{V}{V}{V}{x_1}{0}{y_1}{0}\]
maps an element $f \in
\Multi^K_{\epsfig{file=Images/flat_3.eps,height=3mm}}(V,V,V;V)$ to the
product of vertex operators, $Y(cdot, y_1)Y(\cdot, x_1)\cdot$.  Then the
locality condition reduces to the requirement that $f$ be symmetric (this is
equivalent to the rationality and commutativity of products of Frenkel,
Huang, and Lepowsky, \cite{FHL}).  In fact, the refinement, 
\[\flatthreeleafedlabelled{V}{V}{V}{}{z_1}{z_2}{0} \longrightarrow  
\threeleafedleft{V}{V}{V}{V}{x_1}{0}{y_1}{0}\]
is interpreted as giving a map from $f$ to $Y(Y(\cdot, x_1)\cdot, y_1)\cdot$,
and thus formalises the notion of associativity for a vertex algebra as in
\cite{FHL} or \cite[section 4.6]{kac}.
\end{example}

\begin{proof}[Proof of theorem \ref{t:refinement_works}]
The idea of the proof is that we get a map between pullbacks by showing that
there is a map between the spaces over which we are pulling back.  There is
also a map between the pullback objects, and these maps form a commuting
diagram, inducing a map between pullbacks.  For more of the abstract details
see \cite{cts_valgdetails}.

Explicitly, a tree $p$ can arise as a refinement of a tree $q$ by a sequence
of either replacing nodes of $q$ with edges, or shrinking internal edges of
$q$ down to nodes (see appendix \ref{a:cat_of_trees}).  These two moves can
be interchanged, and can be carried out one move at a time.  Thus we shall
consider them separately.  And because we have defined $\Multi^K_p$ and
$\Multi^K_q$ as an iteration of operators associated to each internal node,
we need only consider refinements involving flat trees.
\begin{description}
\item[Case 1:] We consider the case where $p$ arises as a refinement of $q$
by replacing a node of $q$ with an edge.  Taking $q$ to be the flat $n$
leafed tree, we may either replace the root or a leaf to give $p$.  If we
replace the root, then we are defining a canonical map
\[\generictailA{A}{B}{n} \longrightarrow \genericA{A}{B}{n}\]
Labelling the interior node of the tree on the left hand side, $p_1$, we see
that an arbritrary pullback object in the definition of $\Multi^K_p$, is
\[\Sing_{\sigma,p_1}\Hom_G\bigl(\G{n},\Hom(G,B)\bigr),\]
which is manifestly isomorphic to the pullback object,
$\Sing_{\sigma,\Root}\Hom(\G{n},B)$, of the right hand side.  The same holds
for the space over which we pullback, and so we see that in fact this map is
an isomorphism.  If we instead replace a leaf, we are defining a canonical
map
\[\generictopA{A}{B}{n} \longrightarrow \genericA{A}{B}{n}\]
A similar proof shows that this is also an isomorphism.

\item[Case 2:] We next consider the case where $p$ arises by shrinking an
internal edge of $q$ down to a node.  In this case we are defining a
canonical map
\[\genericA{A}{B}{n} \longrightarrow \genericcompositionB\]
where, without loss of generality, we have chosen to consider refinement of
the first $m$ leaves ($m \leq n$).  Labelling the internal node of the tree
on the right hand side, $q_1$, we see that an arbritrary pullback object is
of the form
\[\Sing_{\sigma,q_1}\Sing_{\delta,\Root_q}\Hom(\G{n-m+1},B),\]
where $\sigma \in A_m$ and $\delta \in A_{n-m+1}$.  We need to show that
there exists a natural map from $\Multi^K_p$ to this object, so considering
an arbritrary element of $\Multi^K_p$, we consider it as an element of the
pullback object $\Sing_{(\sigma,\delta),\Root_p} \Hom(\G{n},B)$, where we
define a permutation, $(\sigma,\delta) \in A_n$:
\begin{eqnarray*}
(1, \ldots, m) &\mapsto & (\sigma(1), \ldots, \sigma(m))\\
(m+1, \ldots, n) &\mapsto & (\delta(x), \delta(m+1),\ldots,
\widehat{\delta(i)}, \ldots, \delta(n)) 
\end{eqnarray*}
where $x$ is the placeholder for the internal edge, and where we exclude
$\widehat{\delta(i)}$ when $i$ is mapped to $x$.  Notice that we have shifted
the domain of $\delta$ so that it acts on $(m,\ldots, n)$, and that the
permutation, $(\sigma,\delta)$, will arise for $n-m+1$ different choices of
$\delta$.

Consider the canonical map we are trying to exhibit:
\[\Sing_{(\sigma,\delta),\Root_p} \Hom(\G{n},B)\longrightarrow
\Sing_{\sigma,q_1}\Sing_{\delta,\Root_q}\Hom(\G{n-m+1},B).\] It will be an
identity on the outer terms corresponding to $\Sing_{\sigma,q_1}$.  So our
problem reduces to showing that there exists a map
\begin{multline}\label{e:refine_problem}
\Hom\biggl(A_{(\sigma,\delta)(m+1)}, K^{\otimes m} \tens \Hom\Bigl(
A_{(\sigma,\delta)(m+2)}, K^{\otimes m+1} \tens \cdots
\Hom(\G{n},B)\Bigr)\biggr) \\
\longrightarrow 
\Hom\Biggl(A_{\delta(x)} \tens A_{\delta(m+1)}, K \tens
\Hom\biggl(A_{\delta(m+1)}, K^{\otimes 2} \tens \Hom\Bigl( A_{\delta(m+2)},
K^{\otimes 3} \tens \cdots \Hom(\G{n-m+1},B)\Bigr)\biggr)\Biggr),
\end{multline}
where $A_x = \G{m}$.  We immediately notice that there are fewer copies of
the singularity $K$ on the right hand side.  This is because when we expand
the singularities, we will have maps of the form:
\[K\tens K \longrightarrow \Hom_G(G,K)\tens \Hom_G(G,K) \longrightarrow
\Hom_G(\G{2},K),\] 
where the second arrow is achieved by multiplication in $K$ (e.g., example
\ref{ex:expand_three}).  A complicated but routine calculation shows that a
succession of maps of this form give the desired map. 
\end{description}
\end{proof}

\section{Relaxed Multi Categories}

So far in this paper, we have introduced the notion of a vertex group.  We
then looked at representations of vertex groups, paying special attention to
representation of the classical vertex group.  We then introduced the notion
of multilinear singular maps, and saw that in order for them to compose
properly, we needed to parameterise them by trees.  Finally we described the
canonical maps between singular maps which arise when we refine one tree to
another.

In this section we shall be concerned with general properties satisfied by
our multilinear singular maps.  The work which we have already done in
defining, composing and refining the singular maps will be enough to make it
clear that we have been working with examples of the following structure.

\begin{dfn}\label{d:relmultcat}
A \bemph{relaxed multilinear category} is an ordinary category of $R$\h
modules for some ring $R$ with the following additional structure:
\begin{enumerate}
\item \bemph{Multi Maps:} For any $n+1$ objects $A_1, \ldots, A_n, B$ and any $n$
leafed tree, $p$, there is a collection of multi maps from $A_1, \ldots, A_n$ to
$B$ denoted $\Multi_p(A_1, \ldots, A_n, B)$.  In particular, if $\bullet$
is the unique tree with no internal vertices, then for all objects $A$,
$\Multi_\bullet(A,A)$ contains the identity morphism.
\item\label{d:relaxed_composition} \bemph{Composition:} For any $(n+1)$\h tuple
of multi maps, $\Multi_{p_i}(A_{i1}, \ldots, A_{im_1}, B_i)$, $\Multi_q(B_1,
\ldots, B_n, C)$, $1 \leq i \leq n$ there is a map that
defines composition: 
\begin{multline*}
\Multi_{p_1}(A_{11},\ldots, A_{1m_1}, B_1)\tens \ldots \tens
\Multi_{p_n}(A_{n1},\ldots, A_{nm_n}, B_n) \\ \tens
\Multi_q(B_1, \ldots, B_n, C) \longrightarrow
\Multi_{q(p_1, \ldots, p_n)}(A_{11},\ldots, A_{nm_n}, C),
\end{multline*}
where $q(p_1, \ldots, p_n)$ is the tree with $\sum m_i$ leaves formed by
gluing the root of each tree $p_i$ to the $i$th external edge.  Composition
of maps of this type is associative.
\item\label{d:relaxed_refinement} \bemph{Refinement:} If $p$ and $q$ are $n$
leaved trees and $p$ is a refinement of $q$, then for every collection of
$n+1$ objects $A_1, \ldots, A_n, B$ there is a map
\[r_{p,q}: \Multi_p(A_1, \ldots, A_n, B) \longrightarrow \Multi_q(A_1, \ldots, A_n, B).\]
\item \bemph{Linearity:} The collections of multi maps, $\Multi_p(A_1,
\ldots, A_n, B)$ are $R$\h modules and composition is multilinear.
\end{enumerate}
\end{dfn}

The category of representations of a vertex group possesses the structure of
a relaxed multilinear category using the $G$\h invariant collections of
singular functions.  Notice that we need to restrict to $G$\h invariant
multi maps in order for the composition condition to hold.  

With this categorical structure defined, the next natural step is to define
an algebra in a relaxed multilinear category.  Recall that we defined an
(associative) algebra in a symmetric tensor category in definition
\ref{d:algebra_defn}.  In a relaxed multilinear category we define an
(associative) algebra as follows:

\begin{dfn}
An \bemph{(associative) algebra} in a relaxed multilinear category, $\cB$,
consists of an object $B \in \cB$ and a collection of maps $\{f_{p}\} =
\set{f_{p} \in \Multi^K_{G,p}(B, \ldots, B, B)}{\text{$p$ is an $n$ leafed
tree},n \in \N}$ .  These maps must satisfy the following axioms:
\begin{enumerate}
\item\label{d:algebra_composition} \bemph{Composition:} If $q(p_1, \ldots,
p_n)$ is the tree formed by gluing the root of each tree $p_i$ to
the $i$th external edge of an $n$ leafed tree, $q$, ($p_i$ possibly empty), then 
\begin{equation}
f_{q(p_1, \ldots, p_n)} = f_q \circ (f_{p_1}, \ldots, f_{p_n}).
\end{equation}
\item \bemph{Unit:} The map $f_\circ:R \rightarrow B$ (where $\circ$ is the
empty tree) defines a unit for the algebra in the sense that for any $n$
leafed tree, $p$, and any $1 \leq k \leq n$,
\[f_p \circ_k f_\circ = f_{p^\prime}\] 
where $\circ_k$ denotes composition at the $k$th leaf of $p$, and  $p^\prime$
is the $n-1$ leaved tree arrived at by removing the $k$th leaf from $p$.
\item \bemph{Refinement:} If $p,q \in \T{n}$ and
$p$ is a refinement of $q$, then $r_{p,q} (f_p) = f_q$ where $r_{p,q}$ is
the refinement map given by the refinement axiom for a relaxed multilinear
category.
\end{enumerate}
\end{dfn}

This is an algebra in the sense that each map $f_p$ defines an ``$n$\h fold
multiplication'' for elements of $B$.  For all $n \in \N$ we denote the
multilinear map associated to the flat tree with $n$ leaves by $f_n$.  Since
composition of multi maps in $\cB$ is associative, the associativity of
$(B,\{f_{p}\})$ is a consequence of the composition axiom.  Considering
$\bullet$, the 1 leafed tree with zero edges, then since $f_p \circ_k
f_{\bullet} = f_p$ and $f_{\bullet} \circ f_p = f_p$, we see that
$f_{\bullet} = 1_B$.  The algebra defined by $(B,\{f_{p}\})$ is said to be
\defn{commutative} if there exists an action of the symmetric group on each
of the multilinear maps in $\{f_{p}\}$.  
\begin{example}\label{ex:last_part_of_vacuum}
Since $f_{\bullet}$ refines to
$f_{\ \epsfig{file=Images/1_flat.eps,height=3mm}}$, we know that for any $b \in
B$, and $g \in G$,
\[f_{\ \epsfig{file=Images/1_flat.eps,height=3mm}}(b)(g) = gb.\]
\end{example}

\begin{note}
This definition of an algebra is just a functor from the opposite of the
category or trees (see appendix \ref{a:cat_of_trees}) to $\cB$ where each
object $p \in \mathcal{T}$ is mapped to an element of $\Multi^K_{G,p}$.
\end{note}

\begin{thm} \label{t:from_alg_to_valg}
Given any commutative algebra, $(B,\{f_p\})$, for the classical vertex group,
it gives rise to a vertex algebra as in definition \ref{d:vetex_algebra}
where $B$ is the space of fields and $\vac = f_0(1)$.  The infinitesimal
translation operator is $T=\D{1}$, and the state-field correspondence is
given by $Y(\cdot,x)\cdot = f_2(\cdot,\cdot)|_{y=0} : B \rightarrow
B\ps{x}[x^{-1}]$.
\end{thm}

\begin{proof}
In order to show that $(B,\{f_p\})$ gives a vertex algebra, we check the
axioms for the vertex algebra.  The vacuum axioms follow naturally from the
discussion of section \ref{ss:vacuum}, together with previous example,
\ref{ex:last_part_of_vacuum}, giving
\[f_{\epsfig{file=Images/2_leaf.eps,height=3mm}}(\vac \tens b) = \sum_{i\geq
0} \frac{T^i}{i!}b y^i = e^{Ty}b.\] 
Translation covariance is automatically satisfied because our functions,
$f_p$, are $G$\h invariant.  And, we saw that the locality condition is
satisfied in section \ref{ss:locality}.
\end{proof}

\begin{thm}
A vertex algebra (as in definition \ref{d:vetex_algebra}) with state space
$V$, defines an algebra in the relaxed multilinear category of
representations of the classical vertex group.
\end{thm}

\begin{proof}
As in the proof for holomorphic vertex algebras (claim
\ref{c:holo_equivalence}), the vacuum defines a map $f_\circ$ for the empty
tree.  We described in example \ref{ex:getting_singmaps} how to construct
$f_{\epsfig{file=Images/2_leaf.eps,height=3mm}}$, and section
\ref{ss:locality} showed how how to construct the singular function
associated to the flat tree with three leaves.  Carrying out a similar
process leads to the construction of singular functions associated to all
flat trees, and the singular functions associated to trees with internal
nodes arise from composition of singular functions associated to flat trees.
Closure under refinement follows by construction. 
\end{proof}

\bigskip
\appendix
\section{Trees}\label{a:cat_of_trees}
Formally, a \defn{tree} is defined to be a connected oriented planar graph
with a finite number of vertices and no circuits.  The empty tree is denoted
by a circle, $\circ$.  All other trees have at least one vertex, and given
any such graph we choose a particular vertex which we call the \defn{root} of
the tree.  For convenience we will always draw the root of the tree at the
bottom of the graph.  The vertices which are joined to exactly one edge
(excluding the root) are called the external vertices or the \defn{leaves} of
the tree.  All other vertices are called internal vertices.  For $n \geq 0$,
we let \T{n} denote the collection of all $n$ leafed trees.  

We say that a tree has \defn{height} $d$, if $d$ is the greatest number of
edges between the root vertex and a leaf.  In \eqref{e:3leafs}, the first
tree has height 1 while the second and third trees have height 2.  The unique
tree with zero internal vertices is a refinement for every tree in \T{n} and
is the only tree with height 1 in that collection.  We call trees of height 1
\defn{flat} tress.  Each collection of trees \T{n} has a subcollection in
which the leaves of each tree are separated from the root vertex by the same
number of edges.  We denote these \defn{the trees of constant height} \TT{n}.

\begin{equation} \label{e:3leafs}
\flatthreeleafed{}{}{}{} \flatthreeleafedleft{}{}{}{} \flatthreeleafeddown
\end{equation}

For any two trees $p,q$ with the same number of leaves, we say that $p$ is a
\defn{refinement} of $q$ if $p$ arises after a succession (possibly zero) of
the following moves:
\begin{itemize}
\item an internal edge is shrunk down to a vertex,
\item a vertex is replaced by an edge.
\end{itemize}
By internal edges, we mean those edges that do not end in a leaf.  Two trees
are considered to be the same exactly when they have the same oriented graph.
In \eqref{e:3leafs} the first tree is a refinement of the second, and the
first and third trees are refinements of one another.

There are a number of (possibly inequivalent) ways of giving the collection
of all trees the structure of a category. For our purposes we take the
definition due to Tom Leinster \cite{leinster_enrichment} which seems to arise most
naturally when dealing with higher dimensional categories.  (For possibly
different definitions see \cite{soibelman} or \cite{kreimer_trees}.)  In this
categorical structure, we define a single morphism from a tree $q$ to a tree
$p$ exactly when $p$ is a refinement of $q$.  So in equation \eqref{e:3leafs}
there exists an arrow from the second tree to the first, and there exists an
arrow in each direction between the first tree and the third.  Refinement is
transitive so morphisms compose and each tree is a refinement of itself so we
have identity morphisms.  For each $n \geq 0$, \T{n} is a category and so
$\mathcal{T}$ is a category with countably many disconnected components.
Under this definition, the subcollection of trees of constant height forms a
full subcategory of the category of trees, and each \TT{n} is a full
subcategory of \T{n}.

The category of trees naturally forms an operad, which is the usual
multi category of composable trees.

\section{Evaluation Maps}\label{a:evaluation}

This appendix is concerned with making clear what types of evaluation maps
arise naturally when dealing with collections of multilinear maps of $G$\h
modules for an arbritrary vertex group $G$.

\begin{example}
If $A$ and $B$ are $G$\h modules, then we know that there exists a natural
evaluation map
\[\Hom(A,B) \tens B \longrightarrow B.\]
This of course holds for ordinary $R$\h linear maps, $\Hom_R(A,B)$, as well as
for $G$\h invariant maps $\Hom_G(A,B)$.  In fact, given $G$ modules, $A_1,
\ldots, A_n, B$ we have partial evaluation maps
\begin{equation}\label{e:partial_eval}
\Hom(A_1 \tens \cdots \tens A_n, B) \tens A_1 \tens A_2 \longrightarrow
\Hom(A_3 \tens \cdots \tens A_n, B).
\end{equation}
\end{example}

The problem with evaluation arises when we consider singular maps.  Recall
that for a $G$\h module, $B$, the collection of singular maps $\Fun(\G{2},B)$
was defined to be $\Hom(\G{2},B)\tens_{f_{1,2}} K$, where the tensor product is
taken over $H^*$.  We are tempted to evaluate this map for some $\G{2}$,
\[\Hom(\G{2},B)\tens_{f_{1,2}} K \tens \G{2} \longrightarrow B \tens K,\]
but no such natural map exists because $B$ does not have the structure of an
$H^*$\h module.  (Recall that the action of $H^*$ on $\Hom(\G{2},B)$ is on
its domain.)  The following example makes the problem more explicit:

\begin{example}
When $G$ is the classical vertex group, $\Fun(\G{2},B) \isom
B\ps{x,y}[(x-y)^{-1}]$, and so an arbritrary element of this collection of
maps is, 
\[(x-y)^{-k}\sum_{i,j \geq 0} b_{i,j} x^i y^j, \]
for some $b_{i,j} \in B$ and some $k \geq 0$.  Naively evaluating this at
$\D{p}\tens\D{q} \in \G{2}$ we have $b_{p,q}(x-y)^{-k} \in B \tens K$.  But
if we now rewrite our power series as
\[(x-y)^{-k-1}\sum_{i,j \geq 0} b_{i,j} (x^{i+1} y^{j} + x^{i} y^{j+1}, \]
we see that this evaluates to $(x-y)^{-k-1}(b_{p-1,q} + b_{p,q-1})$ for $p,q
\geq 1$.  Thus we see that the same power series seems to have evaluated to
two different solutions.  We still might think that the situation could be
reconciled by setting the two evaluation equal to one another, giving
\[b_{p-1,q} + b_{p,q-1} = b_{p,q}(x-y).\]
But this is just a restriction on the form of our power series, it does not
define an action of $H^*$ on $B$ as we might have hoped.  
\end{example}

Now that we see that there are evaluations which we can and can not make, we
would like to know how to evaluate multilinear singular maps,
$\Multi^K_{\epsfig{file=Images/generic.eps,height=3mm}}(A_1, \ldots, A_n;
B)$, for an arbritrary tree, p.  Since these are defined as pullbacks over
$\Hom_{\G{n}}\biggl(A_1\tens \cdots \tens A_n, \Fun_p(\G{n}, B)\biggr)$, we
know that they can be evaluated partially or completely for all $A_1, \ldots,
A_n$ as in equation \eqref{e:partial_eval}.  In fact, the singular maps,
$\Sing_\sigma(A_1, \ldots, A_n; \Hom(\G{n},B))$ of equation
\eqref{e:nary_general_pbk} can also be evaluated, because although there are
many copies of the singularity $K$ among the $A_i$, the tensoring over $H^*$
is with the inner copy of $\G{n}$ (this can be realised as a suitable
quotient), and so we could for example evaluate the following at
$A_{\sigma(2)}, A_{\sigma(3)}$:
\begin{multline}\label{e:sing_eval}
\Hom\biggl(A_{\sigma(1)}\tens
A_{\sigma(2)}, K \tens \Hom\Bigl(A_{\sigma(3)}, K^{\otimes 2}\tens \bigl(\cdots
K^{\otimes n-1} \tens \Hom(\G{n},B)\cdots\bigr)\Bigr)\biggr) \\
\longrightarrow \Hom\biggl(A_{\sigma(1)}, K \tens \Hom\Bigl(A_{\sigma(4)},
K^{\otimes 2}\tens K^{\otimes 3}\tens \bigl(\cdots
K^{\otimes n-1} \tens \Hom(\G{n},B)\cdots\bigr)\Bigr)\biggr).
\end{multline}

\section{Selected Results from Axiomatic Vertex
Algebras}\label{a:axiom_facts} In this appendix we gather a few facts about
ordinary axiomatic vertex algebras which we have used throughout this paper.
For more details see \cite{kac}.  

Let $(V, Y(\cdot,x)\cdot, T, \vac)$ be a vertex algebra.  This first lemma is
a consequence of the vacuum axioms:

\begin{lemma}
For any $a \in V$, $Y(a,x)\vac = \sum_{i\geq 0} \frac{T^i}{i!}a x^i = e^{xT}a.$
\end{lemma}

\begin{proof}
Since $Y(a,x)\vac$ can be evaluated at zero, we know that it has no
singularities.  By the first vacuum axiom we know that $T\vac =0$, so the
translation covariance axiom says that for any $n \geq 0$, $(\partial_x)^n
Y(a,x)\vac = \frac{T^n}{n!}Y(a,x)\vac$.  Letting $Y(a,x)\vac = \sum_{i\geq 0}
a_ix^i$, and evaluating this at zero we have $a_n = \frac{T^n}{n!}a$ and the
lemma is proved.
\end{proof}

The next lemma uses this result, and is a consequence of the locality axiom.

\begin{lemma}[Quasisymmetry]
For any $a,b \in V$ we have $Y(a,x)b = e^{xT}Y(b,-x)a.$
\end{lemma}

\begin{proof}
Recall that the translation covariance axiom says that $Y(a,x)T =
(T-\partial_x)Y(a,x)$.  So we have
\begin{eqnarray*}
Y(b,y)e^{xT}a &=& \sum_{i\geq 0}Y(b,y)\frac{T^i}{i!}a\\
&=& \sum_{p,q\geq 0}(-1)^q\frac{T^p}{p!}\frac{\partial_y^q}{q!}x^{p+q} Y(b,y)a \\
&=& \sum_{p\geq 0}\frac{T^p}{p!}x^p Y(b,y-x)a \\
&=& e^{xT}Y(b,y-x)a.
\end{eqnarray*}
From the locality axiom we know that for some $N \gg 0$ the following holds
\[(x-y)^{N}Y(a,x)Y(b,y)\vac = (x-y)^{N}Y(b,y)Y(a,x)\vac. \]
Combining this with the previous lemma, it just says that
\begin{eqnarray*}
(x-y)^{N}Y(a,x)e^{yT}b &=& (x-y)^{N}Y(b,y)e^{xT}a\\
&=& (x-y)^{N}e^{xT}Y(b,y-x)a.
\end{eqnarray*}
Setting $y$ equal to zero we have $x^{N}Y(a,x)b = x^{N}e^{xT}Y(b,-x)a$, and
since this is an equality of Laurent series, we may divide by $x^{N}$ to give
our result.
\end{proof}

\section{Borcherds' Singular maps parameterised by
trees}\label{ss:bor_sing_maps} 
This appendix provides a review of the singular multilinear maps defined in
\cite{bor}.  In that paper, the definitions which follow were intended only
for trees of constant height (i.e., trees whose root is separated from each
leaf by the same number of edges), which are equivalent to the \defn{sieves}
defined in that paper.  Here we extend the definition of that paper to all
trees.  Our definition for multilinear singular maps uses these singular maps
as a base over which we pullback to get a more specific collection of
multilinear singular maps.

\begin{dfn}
If $p$ is a tree with $n$ leaves and height $d$, then we define the
\defn{space of singular functions of type $p$ from $\G{n}$ to $B$}, written
as $\Fun_p(\G{n},B)$, recursively on height.  If the height of $p$ is one,
then $p$ is the flat tree with zero internal vertices, and we define
$\Fun_{p}(\G{n},B)$ to be $\Fun(\G{n},B)$ as in definition \ref{d:fun}.  If
$d > 1$, then let $q$ be the tree of height $d-1$ with $m$ leaves, formed by
removing the $s$ edges which are separated from the root by $d-1$ edges.  We
define $\Fun_{p}(\G{n},B)$ as follows:
\begin{enumerate}
\item\label{step1}Begin by forming the collection of $\G{m}$\h invariant maps
$\Hom_{\G{m}}(\G{s},\Fun_{q}(\G{m},B))$.  We make $\G{s}$ into a $\G{m}$\h module
by first noticing that each leaf $0 \leq l \leq m$ of $q$ is joined to
$k_{l}$ edges of $p$.  So the $l$th entry of $\G{m}$ acts diagonally on the 
corresponding $k_{l}$ entries of $\G{n}$.  
\item Localise $\Hom_{\G{m}}(\G{n},\Fun_{q}(\G{m},B))$ at all pairs $(i,j)$ of
the $k_{l}$ entries of $\G{n}$ which are joined to each leaf $1 \leq l \leq
m$.  This collection of singular maps is $\Fun_{p}(\G{n},B)$.
\end{enumerate}
\end{dfn}

\begin{note}
This difference between this definition and the one given in \cite{bor} is
that in step \ref{step1} above, we require the maps from $\G{n}$ to
$\Fun_q(\G{m},B)$ to be $G$\h invariant at each leaf of $q$.  We will see in
the following examples that this definition gives the results put forward
in that paper.  It is our belief that this assumption of $\G{m}$\h invariance
is presumed in the definition given in \cite{bor}.  
\end{note}

The definition is made clear by working out a number of examples.

\begin{example}
Let $G$ be a vertex group and $B$ be a $G$\h module.  We begin by looking at
the space of singular functions from $G$ to $B$, parameterised by
non-branching trees.  For the tree, $p_1$ consisting of just one edge, we know
that $\Fun_{p_1}(G,B) = \Fun(G,B) \isom \Hom(G,B)$.  Using the given
recursion definition, if we let $p_i$ be the tree of height $i$, then
$\Fun_{p_2}(G,B) = \Hom_G(G,\Hom(G,B)) \isom \Hom(G,B)$, and so for any $i
\geq 1$, $\Fun_{p_i}(G,B) \isom \Hom(G,B)$.
\end{example}

\begin{example}
Similarly, if we have an $n$\h leafed tree $p$, and if $q$ is the tree formed
from pasting the root of $p$ onto the end of the tree,
\epsfig{file=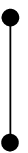,height=3mm}, then $\Fun_{q}(\G{n},B)$ is the
module $\Hom_G(\G{n},\Hom(G,B))$ localised at all pairs $(i,j)$.  But we know
$\Hom_G(\G{n},\Hom(G,B)) \isom \Hom_R(\G{n},B)$, and so $\Fun_{q}(\G{n},B) \isom
\Fun_{p}(\G{n},B)$.
\end{example}

\begin{example}
Let $t_3$ be the unique 3\h leafed tree of height 1.  We saw in example
\ref{ex:G3} that for the classical vertex group, 
\begin{eqnarray*}
\Fun_{t_3}(\G{3},B) &\isom & \Fun(\G{3},B) \\
&\isom & B\ps{x_1,x_2,x_3}[(x_1-x_2)^{-1},(x_1-x_3)^{-1}(x_2-x_3)^{-1}].
\end{eqnarray*}
\noindent If we let the following tree be denoted $l$:
\[\flatthreeleafedleft{}{}{}{}\]
then we have that $\Fun_{l}(\G{3},B)$ is the following
module localised at $(1,2)$:
\[\Hom_{\G{2}}(\G{3},Fun_{t_2}(\G{2},B)) \isom
\Hom_{\G{2}}(\G{3},B\ps{z_1,z_2}[(z_1-z_2)^{-1}]).\] 
Localising, we can identify this module with a subset of the collection of
formal power series 
\begin{equation}\label{e:3leafed_labeled}
B\ps{z_1,z_2}[(z_1-z_2)^{-1}]\ps{y_1,y_2,y_3}
[(y_1-y_2)^{-1}]. 
\end{equation}
The $\G{2}$\h invariance can be realised as the restriction to power series
which satisfy $\partial_{z_1} = \partial_{y_1}+\partial_{y_2}$ and
$\partial_{z_2} = \partial_{y_3}$.  We represent this module pictorially as
the labelled tree,  
\[\flatthreeleafedleftsides{y_1}{y_2}{y_3}{z_1}{z_2}{B}\]
where at each internal node we have equalised over an action of $G$ on the
incoming and outgoing edges.  These differential equations give conditions
which allow us to rewrite this power series module as a power series of only
three variables.  Either of the changes of variables given in equations
\eqref{e:varchange1}-\eqref{e:varchange3} work, and so we may made the
identification
\begin{eqnarray*}
\Fun_l(\G{3},B) &\isom& B\ps{X_1,X_3}[(X_1-X_3)^{-1}]
\ps{(X_1-X_2)}[(X_1-X_2)^{-1}]\\ 
&\isom& B\ps{X_1,X_3}[(X_1-X_3)^{-1}]\ps{X_2}[(X_1-X_2)^{-1}].
\end{eqnarray*}
\end{example}

\bibliography{equivalent}
\nocite{*}
\bibliographystyle{plain}

\end{document}